\documentstyle[12pt,fullpage,emlines,bezier]{article}

\linethickness{0.4pt}

\newcommand{\bi}{\bibitem}
\newcommand{\nb}{\newblock}

\newcommand{\be}[1]{\begin{equation}\label{#1}}
\newcommand{\ee}{\end{equation}}
\newcommand{\xoi}{\{\,x_0^{\pm1},x_1^{\pm1}\,\}}
\newcommand{\la}{\langle\,}
\newcommand{\ra}{\,\rangle}
\newcommand{\prf}{{\bf Proof.}\ }
\newcommand{\pp}{{\cal P}_5}
\newcommand{\iv}{^{-1}}
\newcommand{\bdelta}{\overline{\Delta}}

\newtheorem{thm}{\quad Theorem}
\newtheorem{lm}{\quad Lemma}
\newtheorem{cy}{\quad Corollary}

\title{The Dehn Function of Richard Thompson's Group $F$ is Quadratic}
\author{\vspace{2ex}
V.~S.~Guba\thanks{This research is partially supported by the RFFI
grant 99--01--00894 and the INTAS grant 99--1224.}\\ Vologda State
Pedagogical University,\\ 6 S.~Orlov Street,\\ Vologda\\ Russia\\
160600\\ E-mail: guba{@}uni-vologda.ac.ru}

\begin{document}

\maketitle

\begin{abstract}

We prove that the Dehn function (that is, the smallest
isoperimetric function) of the R.\,Thompson's group $F$ is
quadratic.

\end{abstract}

The Richard Thompson group $F$ can be defined by the following
presentation

\be{xinf} \langle\,x_0,x_1,x_2,\ldots\mid x_j^{x_i}=x_{j+1}\
(i<j)\,\rangle, \ee where $a^b=b^{-1}ab$ by definition. This group
was invented by Richard J. Thompson in 1965 during his work in
$\lambda$-calculus. In 1977--1979 it was rediscovered by Dydak,
Freyd and Heller in their work on homotopy idempotents. Let us
mention a few properties of $F$.

\begin{enumerate}

\item The group $F$ is finitely presented. Namely, it can be given by

\be{x0-1} \langle\,x_0,x_1\mid
x_2^{x_1}=x_3,x_3^{x_1}=x_4\,\rangle, \ee where
$x_{k+1}=x_1^{x_0^k}$ for any $k\ge1$ by definition. So the
defining relations have in fact the following form:
$x_1^{x_0^2}=x_1^{x_0x_1}$, $x_1^{x_0^3}=x_1^{x_0^2x_1}$.

\item Equivalent definition of $F$ can be done in the following way.
Let us consider all strictly increasing continuous piecewise
linear functions from the closed unit interval onto itself. Take
only those of them that are differentiable except at finitely many
dyadic rational numbers and such that all slopes are integer
powers of 2. Obviously, these functions form a group under
composition. This group is exactly $F$.

\item The group $F$ has solvable word and conjugacy problems, it does not
satisfy any non-trivial group law although it has no free
subgroups of rank $>1$. It is known that $F$ is not elementary
amenable. However, the famous problem about amenability of $F$ is
still open.

\end{enumerate}

For the history of the group $F$ see~\cite{CFP}; for the proof of
these and other results on the group $F$ see
also~\cite{BG,BS,Bro,GbS}. \vspace{1ex}

Let us recall the definition of an {\em isoperimetric function\/}
of a group presentation. For our needs it will be sufficient to
give this definition for finite presentations only. However, it
makes sense for any presentation with finite number of generators.
Let $H$ be a group presented by \be{gpfin}
\langle\,a_1,a_2,\ldots,a_s\mid R_1,R_2,\dots,R_t\,\rangle. \ee By
${\cal F}_s$ we denote the free group on $a_1$, $a_2$, \dots,
$a_s$. Let ${\cal N}$ be the normal closure of the defining
relators $R_1$, $R_2$, \dots, $R_t$. By definition, $H\cong{\cal
F}_s/{\cal N}$. For any word $w\in{\cal N}$ we consider the
smallest number $k=k(w)$ such that $w$ is equal in the free group
${\cal F}_s$ to a product of the form $$
U_1^{-1}R_{i_1}^{\pm1}U_1U_2^{-1}R_{i_2}^{\pm 1}U_2\cdots
U_k^{-1}R_{i_k}^{\pm 1}U_k, $$ where $U_i$ are elements of ${\cal
F}_s$. In other words, $k(w)$ is the smallest number of
applications of the defining relations to derive $w$ from the
empty word. Equivalently, using van Kampen's Lemma
(see~\cite{LS}), we may say that $k$ is the smallest number of
cells in a van Kampen diagram over~(\ref{gpfin}) provided its
boundary label is $w$.

A function $p(n)$ is called an {\em isoperimetric function\/} of
presentation~(\ref{gpfin}) whenever $k(w)$ does not exceed $p(n)$
for any word $w\in{\cal N}$ of length $\le n$. The smallest
isoperimetric function of a finite group presentation is called
the {\em Dehn function\/} of this presentation. It is easy to see
that any Dehn function of a group presentation is non-decreasing.

Usually Dehn functions are compared with respect to the following
partial relation $\preceq$ on functions from the set of natural
numbers ${\bf N}$ to itself. By definition, $p\preceq q$ means
that there exists a positive integer constant $C$ such that
$p(n)\le Cq(Cn)+Cn$ for all $n$. This is a pre-order relation on
the set of all non-decreasing functions from ${\bf N}$ to ${\bf
N}$. It induces an equivalence relation $\sim$ on the set of these
functions. Namely, $p\sim q$ if and only if $p\preceq q$ and
$q\preceq p$.

The following result is well known (see~\cite{MO,Ger1,BMS}).
Suppose that $H$ is a finitely presented group and $p$, $q$ are
Dehn functions for two finite presentations of $H$. Then $p$ and
$q$ are equivalent, that is, $p\sim q$. Thus for a finitely
presented group $H$ there exists a unique (modulo the equivalence
relation $\sim$) Dehn function $p(n)$. The group $H$ has solvable
word problem if and only if its Dehn function has a recursive
upper bound.

Gersten proved in~\cite{Ger2} that the Dehn function of $F$ has
exponential upper bound. Also he had a conjecture that the Dehn
function of $F$ is exponential. However, Guba and
Sapir~\cite{GuSa} proved that the Dehn function $\Phi(n)$ of $F$
is strictly subexponential, namely, $\Phi(n)\preceq n^{\log
n}=2^{\log^2n}$. Later Guba~\cite{Gu98} improved this result. He
proved that the Dehn function of $F$ has a polynomial upper bound
$n^5$. Actually it was proved that $\Phi(n)\preceq n^{c}$ for some
constant $c$ between $4$ and $5$.

In the present article we are going to calculate the Dehn function
of $F$ precisely. Here is our main result. \vspace{2ex}

\begin{thm}
\label{main} The Dehn function $\Phi(n)$ of the Richard Thompson
group $F$ is quadratic, that is, $\Phi(n)\sim n^2$.
\end{thm}

Note that any automatic group has quadratic isoperimetric
function~\cite{ECHLPT}. It is still unknown whether $F$ is
automatic.

To prove our main result, we do the following two improvements.
First of all, we improve the estimate for the area of some
``basic" diagrams obtained in~\cite{Gu98}. In that paper it was
$n^{2.746\ldots}$, here we decrease it to $n^2$. There was a lemma
in~\cite{GuSa} that the final estimate for the whole Dehn function
can be obtained by multiplying the previous estimate by $n^2$.
Another important improvement of this paper is to avoid this
multiplier.

Our plan is as follows. In Section~\ref{Stand} we recall some
known facts about string rewriting systems and normal forms in
$F$. Then we define the standard form of each word and reduce the
problem to estimate the Dehn function of $F$ to some fact about
standard diagrams. Namely, we need to show that each word of
length $n$ can be transformed into its standard form in $O(n^2)$
steps. Since the standard form of a word that equals 1 in $F$
freely equals 1, this would imply our main result. The next
Section~\ref{Tri} deals with triangular diagrams. Their boundary
equations have the form $pq=r$, where $p$, $q$ are positive words
with non-decreasing subscripts and $r$ is the standard form of
$pq$. In this Section we reduce the problem for standard diagrams
to the one for triangular diagrams. Section~\ref{Rec} reduces the
problem to rectangular diagrams. Suppose that we conjugate a
positive word $q$ by a letter $x_i$. If all subscripts of $q$ are
greater than $i$ (we say that the conjugation is successful in
this case), then all subscripts increase by 1 after the
conjugation. If we repeat this operation several times, and at
each step the conjugation is successful, then we get a word
obtained from $q$ by increasing all subscripts by $m$, where $m$
is the length of the conjugator $p$. One can draw a diagram of
this conjugation. Such diagrams are called rectangular. We find
two partial cases of rectangular diagrams and formulate two
important statements, Lemma~\ref{hor} and Lemma~\ref{vert}. Now
the problem is reduced to these two cases.

The first of these cases (horizontal diagrams) is studied in
Section~\ref{Ho}. It contains Lemma~\ref{crc}, which is the
crucial point of our proof. Then we go to the other case (vertical
diagrams) in Section~\ref{Ve}. We describe the process to
construct vertical diagrams using auxiliary diagrams with boundary
equation $x_1^{x_0^n}=x_{n+1}$. We need a modification of
Lemma~\ref{crc}. This is Lemma~\ref{intens}. The difference
between it and Lemma~\ref{crc} is that here we need to take areas
of some subdiagrams with coefficients that are called intensities.
After we prove this Lemma, this immediately implies our main
result.

As a corollary, we get that $F$ has linear isodiametric function.
Let us recall the definition of this concept and formulate the
Corollary.

A {\em diameter\/} of a van Kampen diagram is the diameter of its
underlying graph, that is, the greatest distance between its
points. Suppose that for any word $w$ of length $\le n$ that is
equal to 1 over~(\ref{gpfin}), one can find a van Kampen diagram
(over the same presentation) of diameter $\le d(n)$ with boundary
label $w$. The smallest function $d(n)$ with this property is
called the {\em isodiametric function\/} of~(\ref{gpfin}). As in
the case of Dehn functions, there is a similar result that
isodiametric functions of different finite presentations of the
same group are equivalent.

\begin{cy}
\label{isodiam} The isodiametric function of the Richard Thompson
group $F$ is linear.
\end{cy}

This follows from the result of Papasoglu~\cite{Pap}. Note that a
direct proof of the corollary is relatively easy.

At the end of the paper, we discuss possible applications of our
results to some other Thompson-like groups. \vspace{1ex}

The author thanks Mark Sapir for helpful comments.

\section{Standard Forms}
\label{Stand}

We will consider some auxiliary string rewriting system. The facts
about rewriting systems we are going to use are standard. We just
refer to~\cite{DJ90}.

Let $X$ be the alphabet $\{\,x_i,x_i\iv\ (i\ge0)\,\}$. The string
rewriting system $\Sigma$ over $X$ is defined by the following
rewriting rules:

$$ x_i\iv x_i\to x_ix_i\iv,\ \ x_i\iv x_j\to x_{j+1}x_i\iv,\ \
x_j\iv x_i\to x_ix_{j+1}\iv\ \ (i<j), $$

$$ x_jx_i\to x_ix_{j+1},\ \ x_i\iv x_j\iv\to x_{j+1}\iv x_i\iv\ \
(i<j). $$ \vspace{1ex}

Obviously, if $u\to v$ for some words $u$, $v$ over $X$, then
words $u$ and $v$ are equal in $F$, that is, are equal
modulo~(\ref{xinf}). (However, $\Sigma$ is not a rewriting system
for the group $F$.)

It is clear that $\Sigma$ is Noetherian (terminating) because
applying each of its rules decreases words lexicographically,
where $x_0<\cdots<x_n<\cdots<x_n\iv<\cdots<x_0\iv$. The system
$\Sigma$ is also confluent. By Diamond Lemma, it is enough to
check that the system is locally confluent. In our case, this
means that if rules $ab\to cd$, $be\to fg$ belong to $\Sigma$,
where $a,b,c,d,e,f\in X$, then $cde$ and $afg$ can be reduced to
the same word. We shall illustrate this for the case of two rules
$x_i\iv x_j\to x_{j+1}x_i\iv$ and $x_jx_k\to x_kx_{j+1}$, where
$i<j$, $k<j$. We need to show that the words $u=x_{j+1}x_i\iv x_k$
and $v=x_i\iv x_kx_{j+1}$ have a common descendant. There are
three subcases. \vspace{1ex}

a) $i=k$. Here $u\to x_{j+1}x_ix_i\iv\to x_ix_{j+2}x_i\iv$, $v\to
x_ix_i\iv x_{j+1}\to x_ix_{j+2}x_i\iv$.

b) $i<k$. Here $u\to x_{j+1}x_{k+1}x_i\iv\to
x_{k+1}x_{j+2}x_i\iv$, $v\to x_{k+1}x_i\iv x_{j+1}\to
x_{k+1}x_{j+2}x_i\iv$.

c) $i>k$. Here $u\to x_{j+1}x_kx_{i+1}\iv\to
x_kx_{j+2}x_{i+1}\iv$, $v\to x_kx_{i+1}\iv x_{j+1}\to
x_kx_{j+2}x_{i+1}\iv$. \vspace{1ex}

All other cases are quite analogous or they are even easier. We
left the rest to the reader.

Now we can conclude that $\Sigma$ is complete. Therefore, each
word $w$ over $X$ has a unique irreducible form over $\Sigma$. We
call it the {\em standard form\/} of the word $w$.

By a {\em monotone positive word\/} (an MP-{\em word\/}) over $X$
we mean a word of the form $$ x_{i_1}x_{i_2}\cdots x_{i_n}, $$
where $n\ge0$, $i_1\le i_2\le\cdots\le i_n$. \vspace{1ex}

The following easy statement establishes some elementary
properties of the standard form.

\begin{lm}
\label{stf} For any word $w$ over $\xoi$ of length $n$, the
standard form of $w$ is graphically equal to $pq\iv$, where $p$,
$q$ are MP-words, $|pq\iv|=n$ and all subscripts that occur in
$pq\iv$ do not exceed $n$.
\end{lm}

\prf Clearly, all rules of $\Sigma$ preserve the length. If $r$ is
the standard form of $w$, then each positive letter in $r$ is to
the left of each negative letter. Otherwise there is an occurrence
in $r$ of the form $x_i\iv x_i$, $x_i\iv x_j$, or $x_j\iv x_i$,
where $i<j$. This contradicts the fact that $r$ is irreducible. So
$r$ has the form $pq\iv$, where $p$, $q$ are positive words. Since
$r$ is irreducible, it has no occurrences of the form $x_jx_i$,
$x_i\iv x_j\iv$, where $i<j$. Thus $p$ and $q$ are MP-words.

Let us prove by induction on $n$ that maximal subscripts in both
$p$ and $q$ do not exceed $n$. This is obvious if $n=1$. Now
suppose that the fact is true for any word of length $n$. Take any
word $v$ of length $n+1$. It can be presented as $wx_i^{\pm1}$,
$i=0,1$. Let $pq\iv$ be the standard form of $w$, where $p$, $q$
are MP-words with subscripts $\le n$. If $v=wx_0\iv$ then
$p(x_0q)\iv$ is the standard form of $v$. Let $v=wx_0$. The word
$pq\iv$ can be presented as $x_0^sp_1q_1\iv x_0^{-t}$, where
$s,t\ge0$ and $p_1$, $q_1$ do not involve $x_0$. In this case it
is obvious that $pq\iv x_0$ can be reduced to
$x_0^{s+1}p_2q_2^{-1}x_0^{-t}$, where $p_2$ ($q_2$) is obtained
from $p_1$ ($q_1$) by adding 1 to all subscripts. This implies
that all subscripts of the result (which will be the standard form
of $v$) do not exceed $n+1$.

Now let $v=wx_1\iv$. Applying rules of $\Sigma$ to $v$, we move
the last letter to the left whenever possible. Each move may
increase the subscript on this letter by 1 so at the end of the
process the subscript will be at most $n+1$ (the equality holds if
and only if $w=x_0^{-n}$). Also if we move one letter through
another, the letter that stays on the right after the move can
increase its subscript by 1 but this happens only once. So all
subscripts do not exceed $n+1$ when we get the standard form.

The case $v=wx_1$ is analogous.

The proof is complete. \vspace{1ex}

Let us recall some well known facts about normal forms in $F$.
Details can be found in~\cite{CFP}.

Any element of $F$ can be expressed uniquely as a word of the form
\be{NFF} x_{i_1}^{s_1}\cdots x_{i_m}^{s_m}x_{j_n}^{-t_n}\cdots
x_{j_1}^{-t_1}, \ee where $m,n\ge0$, $i_1\le\cdots\le i_m\ne
j_n\ge\cdots\ge j_1$; $s_1,\ldots,s_m,t_1,\ldots,t_n\ge1$. Here it
is also claimed that if both $x_i$ and $x_i\iv$ occur
in~(\ref{NFF}) for some $i\ge0$ then either $x_{i+1}$ or
$x_{i+1}\iv$ also occurs. An expression of the form~(\ref{NFF}) of
an element $g\in F$ is called the {\em normal form\/} of $g$.
(Note that in~\cite{GuSa} it is constructed another useful normal
form for elements of $F$.) In particular, each MP-word is a normal
form. So the word $pq\iv$, where $p$, $q$ are MP-words, may
represent the identity of $F$ if and only if $p$ and $q$ are
graphically equal.

Let \be{p-r} {\cal P}_r=\la x_0,x_1,x_2,\ldots\mid
x_j^{x_i}=x_{j+1} \ (0<j-i\le r)\ra \ee be a group presentation.
It is easy to see that for $r\ge2$ the group presented
by~(\ref{p-r}) is $F$. We will usually construct van Kampen
diagrams over~(\ref{p-r}) to estimate the Dehn function of $F$.
Let us compare ${\cal P}_r$ with finite presentation~(\ref{x0-1}).

By $\psi$ we denote the shift mapping. It takes each letter $x_i$
to $x_{i+1}$ ($i\ge0$). Clearly, $\psi$ induces an embedding of
$F$ into itself. (This notation will be often used throughout the
paper.) We allow to apply $\psi$ and its non-negative powers to
van Kampen diagrams over~(\ref{p-r}). For any diagram $\Delta$
over ${\cal P}_r$, there exists a diagram $\psi^k(\Delta)$ over
${\cal P}_r$ ($k\ge0$) of the same area. (Here we increase all
subscripts of letters by $k$.) The following very elementary fact
shows the relationship between areas of van Kampen diagrams over
$\pp$ and~(\ref{x0-1}). A similar fact is true for any $r\ge2$.

\begin{lm}
\label{areas} Let $w$ be a word over $\xoi$. Suppose that $w$
equals 1 in $F$. Let $\Delta$ be a van Kampen diagram over $\pp$
 of area $N$ with boundary label $w$. Then there exists a van
Kampem diagram over~$(\ref{x0-1})$ with boundary label $w$ and
area $\le13N$.
\end{lm}

\prf We convert each diagram $\Delta$ over $\pp$ into a diagram
$\bdelta$ over~(\ref{x0-1}) in the following way. Each edge
labelled by $x_n$, where $n\ge2$, is replaced by the path labelled
by $x_1^{x_0^{n-1}}$. Each cell of $\Delta$ has boundary equation
of the form $x_j^{x_i}=x_{j+1}$, where $0<j-i\le5$. If $i=0$ then
the boundary label of its image in $\bdelta$ will be freely equal
to 1. Suppose that $i\ge1$. Then the boundary equation of the
image of this cell in $\bdelta$ is $x_k^{x_1}=x_{k+1}$ conjugated
by $x_0^{i-1}$ in the free group, where $k=j-i+1$, $2\le k\le6$.
It is easy to see that the relations $x_4^{x_1}=x_5$,
$x_5^{x_1}=x_6$, $x_6^{x_1}=x_7$ can be derived from the defining
relations of~(\ref{x0-1}) in 5, 9, and 13 steps, respectively.
Recall how to do this. Clearly, any relation of the form
$x_j^{x_i}=x_{j+1}$, where $i\ge1$, $j\in\{\,i+1,i+2\,\}$ is
obtained from one of the defining relations of~(\ref{x0-1}) by
conjugation. So \be{x6} x_6=x_5^{x_3}=x_4^{x_2x_3}=x_3^{x_2^2x_3}
\ee after applying three defining relations. The right-hand side
of this equality is a word of length 7. If we conjugate~(\ref{x6})
by $x_1$ and apply seven more defining relations of~(\ref{x0-1}),
then we obtain $x_6^{x_1}=(x_3^{x_2^2x_3})^{x_1}=x_4^{x_3^2x_4}$.
It remains to note that conjugation by $x_0$ applied to~(\ref{x6})
gives $x_7=x_4^{x_3^2x_4}$ in three steps. So we finally have
$x_6^{x_1}=x_7$ in 13 steps. (The reader can draw the
corresponding van Kampen diagram to see the diagrams for the other
two relations as subdiagrams. They have 5 and 9 cells,
respectively.)

Now after free cancellations of labels inside cells and inserting
diagrams of at most $13$ new cells into each old cell, we get the
diagram $\bdelta$ over~(\ref{x0-1}) with at most $13N$ cells,
where $N$ is the area of $\Delta$. The boundary label remains the
same since $w$ is a word over $\xoi$.

The proof is complete. \vspace{1ex}

Let $w$ be a word over $\xoi$ of length $n$. Let $pq\iv$ be its
standard form. Our aim is to construct a diagram over $\pp$ of
area $O(n^2)$ with boundary equation $w=pq\iv$. This will
immediately imply that the Dehn function $\Phi(n)$ of $F$ is
quadratic. Indeed, if $w$ equals 1 in $F$, then $p$ and $q$ must
be graphically equal. So this gives a diagram over $\pp$ of area
$O(n^2)$ with boundary label $w$. By Lemma~\ref{areas}, there
exists a diagram over~(\ref{x0-1}) of area $O(n^2)$ with boundary
label $w$. Hence $\Phi(n)\preceq n^2$. Since $F$ is not
hyperbolic, $n^2\preceq\Phi(n)$ (this can be easily shown
directly). As a result, $\Phi(n)\sim n^2$, that is, the Dehn
function of $F$ is quadratic. So we proved the following

\begin{lm}
\label{Dehn} Suppose that there exists a positive constant $C$
such that for any word $w$ over $\xoi$, there exists a van Kampen
diagram over $\pp$ of area $\le Cn^2$ with boundary equation
$w=pq^{-1}$, where $pq\iv$ is the standard form of $w$. Then the
Dehn function of $F$ is quadratic.
\end{lm}

\section{Triangular Diagrams}
\label{Tri}

Lemma~\ref{Dehn} says that the problem about the Dehn function of
$F$ can be reduced to standard diagrams, that is, to the diagrams
with boundary equation $w=pq\iv$, where $pq\iv$ is the standard
form of $w$. Now we are going to define some other class of
diagrams and reduce the problem to it. (For any class of diagrams
we deal with, the problem is to show that diagrams over $\pp$ from
this class satisfy quadratic isoperimetric inequality.)

Let $p$, $q$ be MP-words and let $r$ be the standard form of $pq$.
Clearly, $r$ is also an MP-word. We are going to construct a
diagram with boundary equation $pq=r$ over $\pp$ and estimate its
area. In Lemma~\ref{trian} we reduce our problem to this kind of
diagrams that we call {\em triangular diagrams}.

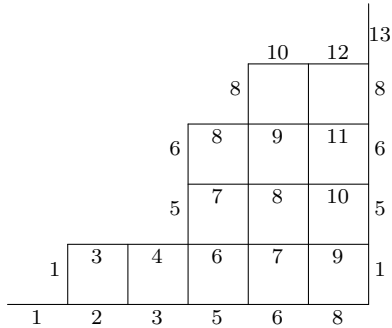
\begin{center}
\unitlength=0.80mm \special{em:linewidth 0.4pt}
\linethickness{0.4pt}
\begin{picture}(65.00,55.00)
\emline{1.00}{5.00}{1}{61.00}{5.00}{2}
\emline{61.00}{5.00}{3}{61.00}{55.00}{4}
\emline{11.00}{5.00}{5}{11.00}{15.00}{6}
\emline{11.00}{15.00}{7}{61.00}{15.00}{8}
\emline{21.00}{5.00}{9}{21.00}{15.00}{10}
\emline{31.00}{5.00}{11}{31.00}{35.00}{12}
\emline{31.00}{35.00}{13}{61.00}{35.00}{14}
\emline{41.00}{5.00}{15}{41.00}{45.00}{16}
\emline{41.00}{45.00}{17}{61.00}{45.00}{18}
\emline{51.00}{5.00}{19}{51.00}{45.00}{20}
\emline{31.00}{25.00}{21}{61.00}{25.00}{22}
\put(6.00,2.00){\makebox(0,0)[cc]{$^1$}}
\put(16.00,2.00){\makebox(0,0)[cc]{$^2$}}
\put(26.00,2.00){\makebox(0,0)[cc]{$^3$}}
\put(36.00,2.00){\makebox(0,0)[cc]{$^5$}}
\put(46.00,2.00){\makebox(0,0)[cc]{$^6$}}
\put(56.00,2.00){\makebox(0,0)[cc]{$^8$}}
\put(46.00,46.00){\makebox(0,0)[cc]{$^{10}$}}
\put(56.00,46.00){\makebox(0,0)[cc]{$^{12}$}}
\put(63.00,10.00){\makebox(0,0)[cc]{$^1$}}
\put(63.00,20.00){\makebox(0,0)[cc]{$^5$}}
\put(63.00,30.00){\makebox(0,0)[cc]{$^6$}}
\put(63.00,40.00){\makebox(0,0)[cc]{$^8$}}
\put(63.00,49.00){\makebox(0,0)[cc]{$^{13}$}}
\put(9.00,10.00){\makebox(0,0)[cc]{$^1$}}
\put(29.00,30.00){\makebox(0,0)[cc]{$^6$}}
\put(16.00,12.00){\makebox(0,0)[cc]{$^3$}}
\put(26.00,12.00){\makebox(0,0)[cc]{$^4$}}
\put(36.00,12.00){\makebox(0,0)[cc]{$^6$}}
\put(46.00,12.00){\makebox(0,0)[cc]{$^7$}}
\put(56.00,12.00){\makebox(0,0)[cc]{$^9$}}
\put(29.00,20.00){\makebox(0,0)[cc]{$^5$}}
\put(39.00,40.00){\makebox(0,0)[cc]{$^8$}}
\put(36.00,22.00){\makebox(0,0)[cc]{$^7$}}
\put(46.00,22.00){\makebox(0,0)[cc]{$^8$}}
\put(56.00,22.00){\makebox(0,0)[cc]{$^{10}$}}
\put(36.00,32.00){\makebox(0,0)[cc]{$^8$}}
\put(46.00,32.00){\makebox(0,0)[cc]{$^9$}}
\put(56.00,32.00){\makebox(0,0)[cc]{$^{11}$}}
\end{picture}
\end{center}

It is useful to describe the structure of triangular diagrams
over~(\ref{xinf}). Let $p$, $q$, $r$ be as above, $q=x_{k_1}\cdots
x_{k_m}$. We define a sequence of words $v_0$, $u_1$, $v_1$,
\dots, $u_m$, $v_m$ by induction. By definition, $v_0=p$. For any
$1\le i\le m$, let $u_i$ be the longest suffix of $v_{i-1}$ that
starts from a letter with a subscript greater than $k_i$. (Here
$u_i$ is empty if no such a letter exists.) By $v_i$ we mean the
word $\psi(u_i)$. Clearly, there exists a rectangular diagram
over~(\ref{xinf}) with $u_i$ at the bottom, $v_i$ at the top,
$x_{k_i}$ on the left and on the right, oriented upwards. Gluing
all these diagrams in a natural way, we get the diagram with
boundary equation $pq=r$, where $p$ is on the bottom, $q$ is on
the right, $r$ is the rest of the boundary. The above picture
illustrates the process for $p=x_1x_2x_3x_5x_6x_8$,
$q=x_1x_5x_6x_8x_{13}$,
$r=x_1^2x_3x_4x_5x_6x_8^2x_{10}x_{12}x_{13}$. (We show only
subscripts of letters; all edges are oriented either from the left
to the right, or upwards.)

Let $w$ be a word over $X$. By $\mu(w)$ we denote the maximal $i$
such that $x_i$ or $x_i\iv$ occurs in $w$. If $w$ is empty, then
$\mu(w)=0$ by definition. It is not hard to see from the
construction of the triangular diagram over~(\ref{xinf}) that
$\mu(r)=\max(\mu(q),\mu(p)+|q|)$. We will use this formula in the
proof of the next lemma, which says that if triangular diagrams
satisfy quadratic isoperimetric inequality, then so do standard
diagrams.

\begin{lm}
\label{trian} Suppose that there exists a constant $C_0>0$ such
that for any MP-words $p$, $q$ there exists a van Kampen diagram
over $\pp$ of area $\le C_0n^2$ with boundary equation $pq=r$,
where $r$ is the standard form of $pq$, provided $|r|\le n$ and
$\mu(r)\le n$. Then there exists a positive constant $C$ such that
for any word $w$ over $\xoi$, there exists a van Kampen diagram
over $\pp$ of area $\le Cn^2$ with boundary equation $w=pq^{-1}$,
where $pq\iv$ is the standard form of $w$.
\end{lm}

\prf Let $|w|=n$. If $n=1$, then $w$ coincides with its standard
form, so there is nothing to prove. We proceed by induction on
$n$. Let $n>1$. We can decompose $w$ into two factors of almost
equal length, namely, $w=w_1w_2$, where
$|w_1|=\left[\frac{n}2\right]$, $|w_2|=\left[\frac{n+1}2\right]$.
Clearly, $|w_1|<|w|$ and $|w_2|<|w|$. Thus we can apply the
inductive assumption and find diagrams $\Delta_i$ ($i=1,2$) over
$\pp$ with boundary equations $w_i=p_iq_i\iv$, where the area of
$\Delta_i$ does not exceed $C|w_i|^2$. It is clear that $$
|w_1|^2+|w_2|^2=\left[\frac{n}2\right]^2+\left[\frac{n+1}2\right]^2=
\left[\frac{n^2+1}2\right]\le\frac{n^2+1}2. $$ Now take the word
$q_1\iv p_2$ and find its standard form $p_3q_3\iv$. Equality
$q_1\iv p_2=p_3q_3\iv$ holds in $F$ so $q_1p_3=p_2q_3$ in $F$. It
is easy to see that the words $q_1p_3$ and $p_2q_3$ have the same
standard form $r$. Indeed, standard forms of both words are
positive. They are equal in $F$ so they must be graphically equal.
Now let $p$ be the standard form of $p_1p_3$ and let $q$ be the
standard form of $q_2q_3$. We can construct 4 triangular diagrams
over $\pp$ that have boundary equations $q_1p_3=r$, $p_2q_3=r$,
$p_1p_3=p$, $q_2q_3=q$. All rewriting rules of $\Sigma$ preserve
the number of positive and negative letters. So $|p_2|=|p_3|$,
$|q_1|=|q_3|$. Thus
$|pq\iv|=|p_1|+|p_3|+|q_2|+|q_3|=|p_1|+|p_2|+|q_1|+|q_2|=|w_1|+|w_2|=
|w|=n$. This also implies $|r|=|q_1|+|p_3|\le n$. It remains to
estimate the maximal subscripts of $p$, $q$, $r$.

From Lemma~\ref{stf} we know that $\mu(p_1)\le|w_1|$. Also we know
that $|p_3|=|p_2|\le|w_2|$. Therefore, $\mu(p_1)+|p_3|\le n$.
Recall that $p_3q_3\iv$ was obtained as the standard form of
$q_1\iv p_2$. By induction on $|q_1|$, it is very easy to show
that $\mu(p_3)$ does not exceed $\mu(p_2)+|q_1|$. This implies
$\mu(p_3)\le n$. By the formula obtained before this Lemma, we
have $\mu(p)\le n$. Similarly, $\mu(q)\le n$. The same
inequalities also show that $\mu(r)\le n$.

By the condition of our Lemma, we can assume that each of the 4
triangular diagrams over $\pp$ has area at most $C_0n^2$. So the
diagram

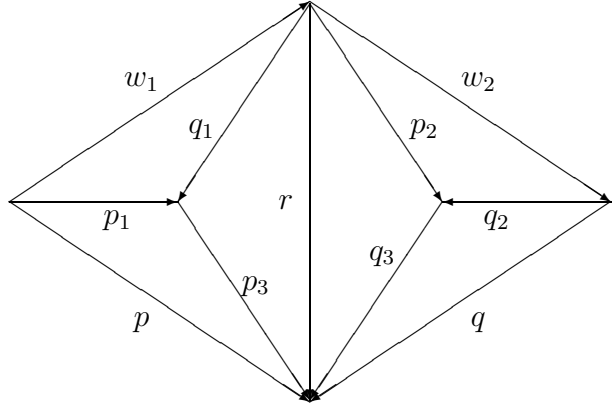
\begin{center}
\unitlength=0.8mm \special{em:linewidth 0.4pt}
\linethickness{0.4pt}
\begin{picture}(101.00,68.33)
\put(1.00,35.00){\vector(3,2){50.00}}
\put(51.00,68.33){\vector(3,-2){50.00}}
\put(1.00,35.00){\vector(3,-2){50.00}}
\put(101.00,35.00){\vector(-3,-2){50.00}}
\put(51.00,68.00){\vector(0,-1){66.00}}
\put(51.00,68.00){\vector(-2,-3){22.00}}
\put(29.00,35.00){\vector(2,-3){22.00}}
\put(51.00,68.00){\vector(2,-3){22.00}}
\put(73.00,35.00){\vector(-2,-3){22.00}}
\put(1.00,35.00){\vector(1,0){28.00}}
\put(101.00,35.00){\vector(-1,0){28.00}}
\put(23.00,55.00){\makebox(0,0)[cc]{$w_1$}}
\put(79.00,55.00){\makebox(0,0)[cc]{$w_2$}}
\put(23.00,15.00){\makebox(0,0)[cc]{$p$}}
\put(79.00,15.00){\makebox(0,0)[cc]{$q$}}
\put(47.00,35.00){\makebox(0,0)[cc]{$r$}}
\put(19.00,32.00){\makebox(0,0)[cc]{$p_1$}}
\put(82.00,32.00){\makebox(0,0)[cc]{$q_2$}}
\put(33.00,47.00){\makebox(0,0)[cc]{$q_1$}}
\put(70.00,47.00){\makebox(0,0)[cc]{$p_2$}}
\put(63.00,26.00){\makebox(0,0)[cc]{$q_3$}}
\put(42.00,21.00){\makebox(0,0)[cc]{$p_3$}}
\end{picture}
\end{center}

\noindent over $\pp$ with boundary equation $w=pq\iv$ has area at
most $C\,\frac{n^2+1}2+4C_0n^2\le Cn^2$ if $C\ge11C_0$ (recall
that $n\ge2$). It remains to mention that $pq\iv$ is the
descendant of $w$ with respect to $\Sigma$ so it is the standard
form of $w$. This completes the proof.

\section{Rectangular Diagrams}
\label{Rec}

Lemmas~\ref{trian} and \ref{Dehn} imply that if triangular
diagrams satisfy quadratic isoperimetric inequality, then the main
result follows. Now we are going to introduce so-called {\em
rectangular diagrams\/} and reduce the problem to them. Let
$p=x_{i_1}\cdots x_{i_m}$, $q=x_{j_1}\cdots x_{j_k}$ be two
MP-words. We are going to conjugate $q$ by $p$. Suppose that
$j_1>i_1$. Then the word $q^{x_{i_1}}$ is equal to
$x_{j_1+1}\cdots x_{j_k+1}$ in $F$. Now suppose that $j_1+1>i_2$.
Then we can conjugate $q^{x_{i_1}}$ by $x_{i_2}$ to get the word
$q=x_{j_1+2}\cdots x_{j_k+2}$, and so on. Thus, if for any $s$
($0\le s<m$) we have $j_1+s>i_{s+1}$, then equality
$q^p=\psi^m(q)$ holds in $F$. One can easily draw a diagram
over~(\ref{xinf}) for this equality (it will be a rectangle with
$km$ cells). In this case we say that conjugation of $q$ by $p$ is
{\em successful\/}.

Suppose now that $p$, $q$ are MP-words such that conjugation of
$q$ by $p$ is successful and let $p=x_{i_1}\cdots x_{i_m}$,
$q=x_{j_1}\cdots x_{j_k}$. Assume that $k$, $m$, $j_k-i_1$ do not
exceed $n$. We say that {\em rectangular diagrams satisfy
quadratic isoperimetric inequality\/} if for any MP-words with the
above conditions there exists a diagram over $\pp$ with boundary
equation $q^p=\psi^m(q)$ of area at most $Ln^2$, where $L>0$ is a
constant. The next Lemma reduces the problem to the case of
rectangular diagrams.

\begin{lm}
\label{rectan} If rectangular diagrams satisfy quadratic
isoperimetric inequality, then so do triangular diagrams.
\end{lm}

\prf For any MP-word $w$ we define its {\em size\/} $||w||$ as the
sum of its length and the difference between its greatest and
smallest subscript. (If
$r=x_1^2x_3x_4x_5x_6x_8^2x_{10}x_{12}x_{13}$, then we have
$||r||=|r|+13-1=23$.) The size of the empty word is zero by
definition.

Note that if an MP-word $w$ is a product of three factors
$w=w_1w_2w_3$, then obviously $||w_1||+||w_3||\le||w||$.

Let $p$, $q$ be MP-words, let $r$ be the standard form of $pq$,
and let $|r|\le n$, $\mu(r)\le n$. Obviously, $||r||\le2n$. It
suffices to show that one can find a van Kampen diagram over $\pp$
of area $O(||r||^2)$ with boundary equation $pq=r$.

We can construct a triangular diagram $T$ with boundary equation
$pq=r$ over~(\ref{xinf}) as we did before. (It is useful to look
at the picture before Lemma~\ref{trian}.) Each letter in $r$ can
be naturally called {\em horizontal\/} or {\em vertical}. If $T$
has at least one cell, then $r$ has an occurrence of the form
$ab$, where $a$ is vertical and $b$ is horizontal. (In the
picture, $r$ has 3 occurrences of this form: $x_1x_3$, $x_6x_8$,
and $x_8x_{10}$.) For each occurrence of this form one can
naturally find a rectangular subdiagram in $T$. (Say, for $x_6x_8$
this rectangular subdiagram will have boundary equation
$(x_5x_6x_8)^{x_1x_5x_6}=x_8x_9x_{11}$.) Denote this rectangular
subdiagram in $T$ by $R$. Taking $R$ off $T$, we get two
triangular diagrams, $T_1$ and $T_2$, where $T_1$ is bounded by
$p$, $r$ and the left side of $R$ whereas $T_2$ is bounded by $q$,
$r$ and the top side of $R$. By $r_1$ and $r_2$ we define the
parts of $r$ that belong to $T_1$, $T_2$, respectively. (In our
example, $r_1=x_1^2x_3x_4$, $r_2=x_8x_{10}x_{12}x_{13}$.) Note
that one of these triangular subdiagrams or both may be empty. If
$T_1$ ($T_2$) is not empty, then $r_1$ ($r_2$) ends (starts) with
a horizontal (vertical) letter. The subdiagrams $R$, $T_1$, $T_2$
and the words $r_1$, $r_2$ depend of the choice of the subword
$ab$. We may sometimes write $r_1(ab)$, $r_1(cd)$ etc if we want
to change the subword.

Estimating areas, we may assume without loss of generality that
$r$ starts with a vertical edge and ends with a horizontal edge.
(In our example we just take the first $x_1$ and $x_{13}$ off
$r$.) This may only decrease $||r||$. Our aim is to find a subword
$ab$ of $r$ in such a way that $||r_1||,||r_2||\le||r||/2$. Let us
show this is always possible. Consider all occurrences of the form
$ab$, where $a$ is vertical, $b$ is horizontal and
$||r_1||\le||r||/2$. Subwords with these properties always exist
because $r_1$ is empty for the leftmost one of them. Now let $ab$
be the rightmost subword of $r$ with the properties listed above.
We claim that $||r_2||\le||r||/2$.

Assume the contrary. Obviously, $r_2$ must have a subword of the
form $cd$, where $c$, $d$ are vertical and horizontal letters,
respectively. (Otherwise $r_2$ is empty.) Choose the leftmost of
these subwords in $r_2$. We have a decomposition $r=r_1r'r_2$,
where $ab$ is contained in $r'$. Elementary properties of the size
imply $||r_1r'||\le||r||-||r_2||<||r||/2$. Let us consider the
rectangular subdiagram $R(cd)$. It is clear that
$r_1(cd)=r_1(ab)r'$. This contradicts the fact that $ab$ was
chosen rightmost.

By induction, we may assume that diagrams $T_1$, $T_2$ can be
filled by cells of $\pp$ in such a way that the area of each of
them will not exceed $K(||r||/2)^2$, where $K>0$ is a constant.
The sum of the horizontal side of $R$ and the vertical side of $R$
does not exceed $||r||$. Also the difference between the greatest
and the smallest subscript that occur in $R$ does not exceed the
difference between the greatest and the smallest subscript of $r$.
Thus the condition of our Lemma allows to fill $R$ by cells of
$\pp$ in such a way that the area of the corresponding rectangular
diagram will not exceed $L||r||^2$. So we get a diagram over $\pp$
with boundary equation $pq=r$. (It consists of $T_1$, $T_2$, and
$R$.) Its area does not exceed $$ 2K\frac{||r||^2}4+L||r||^2\le
K||r||^2 $$ if $K\ge2L$. This completes the proof. \vspace{2ex}

Now it remains to show that rectangular diagrams satisfy quadratic
isoperimetric inequality. Let us consider two partial cases of
this problem. \vspace{1ex}

Case 1 (horizontal diagrams). Let $q$ be an MP-word and let
$p=x_i$ consist of one letter. Assume that the smallest subscript
of $q$ is greater than $i$. Then $q^p=\psi(q)$
modulo~(\ref{xinf}). Suppose that $|q|\le n$ and let $j-i\le n$,
where $j$ is the greatest subscript that occurs in $q$. We need to
construct a diagram over $\pp$ with boundary equation
$q^{x_i}=\psi(q)$ of area $O(n^2)$. Since all subscripts that
occur in $q$ exceed $i$, the word $q$ equals in $F$ to a word of
the form $v(x_{i+1},x_{i+2})$. Suppose that we found a word $v$
such that there exists a van Kampen diagram $\Delta$ over $\pp$ of
area $O(n^2)$ with boundary equation $q=v(x_{i+1},x_{i+2})$. We
can assume that the path labelled by $v$ is simple (otherwise we
can make $v$ shorter). It is easy to see that $v$ cannot be very
long. Indeed, the perimeter of each cell equals 4. So the length
of $v$ cannot be bigger than the number of cells in $\Delta$
multiplied by 4. Thus $|v|=O(n^2)$. (In fact we will find the word
$v$ of linear length.) Conjugating $x_{i+1}$ or $x_{i+2}$ by $x_i$
corresponds to a defining relation of $\pp$. Thus we can construct
a diagram $\Gamma$ over $\pp$ of area $|v|$ with boundary equation
$v(x_{i+1},x_{i+2})^{x_i}=v(x_{i+2},x_{i+3})$. Let us consider the
diagram $\psi(\Delta)$. It has boundary equation
$\psi(q)=v(x_{i+2},x_{i+3})$ and the same area as $\Delta$. It
remains to glue together $\Delta$, $\Gamma$ and the mirror copy of
$\psi(\Delta)$ to obtain the desired ``horizontal" diagram.

To find a word $v(x_{i+1},x_{i+2})$ with the desired properties,
one can first decrease all subscripts of $q$ by $i+1$ (this is
possible since all of them exceed $i$), then express the result as
a word of the form $v(x_0,x_1)$ and then increase all subscripts
by $i+1$. This leads to the following statement we are going to
prove later.

\begin{lm}
\label{hor} For any MP-word $q$ such that $|q|\le n$, $\mu(q)\le
n$, there exists a word $v(x_0,x_1)$ of length $O(n)$ and a van
Kampen diagram over $\pp$ of area $O(n^2)$ with boundary equation
$q=v(x_0,x_1)$.
\end{lm}

Case 2 (vertical diagrams). Now let $q=x_j$ consist of one letter
and let $p=x_{i_1}x_{i_2}\cdots x_{i_m}$ be an MP-word. Suppose
that conjugation of $q$ by $p$ is successful, that is, $j>i_1$,
$j+1>i_2$, $j+2>i_3$, \dots, $j+m-1>i_m$. In this case
$x_j^p=x_{j+m}$ in $F$. We need to find a van Kampen diagram
$\Gamma$ over $\pp$ of area $O(n^2)$ for this boundary equation
provided $j+m\le n$. In fact we need more. Suppose that we have
$l$th power of $x_j$ instead of $x_j$. If we just glue $l$ copies
of $\Gamma$ together, the area becomes $l$ times bigger. To save
the area, we need to mention that $\Gamma$ should have symmetric
structure with respect to a vertical axis. Namely, we want
$\Gamma$ to consist of three parts, the {\em left}, the {\em
central}, and the {\em right}, such that the central part has
$O(n)$ cells and the other two parts are mirror copies of each
other. In this case, if we take $l$th power, this will increase
the area by $O(nl)$ only. To formulate the main fact about
vertical diagrams, we need a definition. \vspace{1ex}

A word $t=x_{j_1}^{d_1}\cdots x_{j_h}^{d_h}$ is called a {\em
smooth\/} word of rank $j$, where $d_i=\pm1$ for all $1\le i\le
h$, whenever the following conditions hold: a) if $d_i=1$, then
$j+d_1+\cdots+d_{i-1}$ belongs to $\{\,j_i+1,j_i+2\,\}$; b) if
$d_i=-1$, then $j+d_1+\cdots+d_{i-1}$ belongs to
$\{\,j_i+2,j_i+3\,\}$. The number $m=d_1+\cdots+d_h$ is called the
{\em height\/} of $t$.

To clarify this definition, let us show how to construct a van
Kampen diagram over $\pp$ with boundary equation $x_j^t=x_{j+m}$.
First of all, let $s_0=j$, $s_1=j+d_1$, \dots,
$s_i=j+d_1+\cdots+d_i$, \dots, $s_h=j+d_1+\cdots+d_h=j+m$. For
each $1\le i\le h$ such that $d_i=1$ we can take a cell over $\pp$
with boundary equation $x_{s_{i-1}}^{x_{j_i}}= x_{s_{i-1}+1}$. If
$d_i=-1$, then we can take the cell with boundary equation
$x_{s_{i-1}}^{x_{j_i}\iv}=x_{s_{i-1}-1}$. In any case, conjugating
$x_{s_{i-1}}$ by $x_{j_i}^{d_i}$ gives $x_{s_i}$. So if we glue
all these $h$ cells together in a natural way, then we get a
diagram over $\pp$ with boundary equation $x_j^t=x_{j+m}$.

\begin{lm}
\label{vert} Let $q=x_j$, $p=x_{i_1}x_{i_2}\cdots x_{i_m}$.
Suppose that $j+m\le n$, $i_1\le i_2\le\cdots\le i_m$ and for each
$0\le s<m$ inequalities $j+s>i_{s+1}$ hold, that is, conjugation
of $q$ by $p$ is successful. Then there exists a smooth word $t$
of rank $j$, height $m$ and a van Kampen diagram over $\pp$ with
boundary equation $p=t$ of area $O(n^2)$. In addition, $|t|=O(n)$.
\end{lm}

First of all let us show that Lemmas~\ref{hor} and~\ref{vert}
imply that rectangular diagrams satisfy quadratic isoperimetric
inequality. According to Lemmas~\ref{rectan}, \ref{trian} and
\ref{Dehn}, this will imply our main result about the Dehn
function of $F$.

Suppose that Lemmas~\ref{hor} and~\ref{vert} hold. Take any
MP-words $p$ and $q$ such that conjugation of $q$ by $p$ is
successful. Let $|p|$, $|q|$ do not exceed $n$ and let the
difference between the greatest subscript of $q$ and the smallest
subscript of $p$ also does not exceed $n$. Our aim is to construct
a van Kampen diagram over $\pp$ of area $O(n^2)$ with boundary
equation $q^p=\psi^m(q)$, where $m=|p|$.

By $j$ we denote the smallest subscript that occurs in $q$, that
is, $q$ starts with $x_j$. We can thus apply $\psi^{-j}$ to $q$.
Obviously, $|\psi^{-j}(q)|=|q|\le n$. The greatest subscript of
$\psi^{-j}(q)$ equals $\mu(q)-j<\mu(q)-i\le n$, where $i$ is the
smallest subscript of $p$. Then we can apply Lemma~\ref{hor} and
find a word $v$ of length $O(n)$ such that some van Kampen diagram
over $\pp$ has boundary equation $\psi^{-j}(q)=v(x_0,x_1)$ and the
area of this diagram is $O(n^2)$. Adding $j$ to all subscripts, we
get a diagram $\Delta$ over $\pp$ of the same area with boundary
equation $q=v(x_j,x_{j+1})$.

By definition, conjugation of $x_j$ by $p$ is also successful. By
Lemma~\ref{vert}, there exists a smooth word $t$ of rank $j$,
height $m$ and a van Kampen diagram $\Gamma$ over $\pp$ with
boundary equation $p=t$ of area $O(n^2)$. Besides, $|t|=O(n)$. By
definition of a smooth word, one can form a diagram $\Xi_0$ over
$\pp$ with boundary equation $x_j^t=x_{j+m}$ that has $|t|$ cells.
Let us draw this diagram vertically, that is, $t$ will be a
vertical path. Take all horizontal edges and increase their
subscripts by 1. This gives some diagram $\Xi_1$. Note that all
defining relations used in $\Xi_0$ had difference 1 or 2. (By a
{\em difference\/} of a defining relation of the form
$x_j^{x_i}=x_{j+1}$ we mean the number $j-i$.) Thus in $\Xi_1$
these differences will be 2 or 3. Presentation $\pp$ allows any
difference from 1 to 5 so $\Xi_1$ will be a diagram over $\pp$.
The number of its cells is also $|t|$.

Now for each $d\in\{\,0,1\,\}$ we have a diagram $\Xi_d$ of $|t|$
cells that has boundary equation $x_{j+d}^t=x_{j+m+d}$. Thus for
any word $v(x_j,x_{j+1})$ there exists a diagram of $|t|\cdot|v|$
cells with boundary equation
$v(x_j,x_{j+1})^t=v(x_{j+m},x_{j+m+1})$. Denote it by $\Xi$. The
number of cells in it is $O(n^2)$ since both $|t|$ and $|v|$ are
$O(n)$.

Let us draw $\Xi$ in such a way that paths labelled by $t$ go
upwards. Take two copies of $\Gamma$ and attach them to $\Xi$
along $t$ on the left and on the right. The boundary equation of
the result will have the form
$v(x_j,x_{j+1})^p=v(x_{j+m},x_{j+m+1})$. Now let us attach
$\Delta$ on the bottom along the word $v(x_j,x_{j+1})$ (recall
that its boundary equation was $q=v(x_j,x_{j+1})$) and then attach
the mirror copy of $\psi^m(\Delta)$ on the top along the word
$v(x_{j+m},x_{j+m+1})$. The final result of all these operations
will have boundary equation $q^p=\psi^m(q)$, as desired. The area
will be $O(n^2)$ since it consists of 5 parts of area $O(n^2)$.
This shows that rectangular diagrams satisfy quadratic
isoperimetric inequality.

So we reduced everything to Lemmas~\ref{hor} and \ref{vert}.

\section{Horizontal Diagrams}
\label{Ho}

We begin to prove Lemma~\ref{hor}. Take an arbitrary MP-word $q$
such that $|q|\le N$, $\mu(q)\le N$. We can write
$q=x_0^kx_{i_1}^{d_1}\cdots x_{i_m}^{d_m}$, where $k\ge0$,
$d_1,\ldots,d_m\ge1$, $0<i_1<\cdots<i_m$. Recall that $x_j$ equals
$x_1^{x_0^{j-1}}$ modulo~(\ref{xinf}) for any $j\ge1$. If we
replace each $x_j$ ($j\ge1$) by $x_1^{x_0^{j-1}}$ and cancel the
result in the free group, then we get the following word $Q$ over
$\xoi$: $$ Q=x_0^{k-i_1+1}x_1^{d_1}x_0^{-(i_2-i_1)}x_1^{d_2}\cdots
x_0^{-(i_m-i_{m-1})}x_1^{d_m}x_0^{i_m-1}. $$

The length of $Q$ equals $$
|Q|=|k-i_1+1|+(i_2-i_1)+\cdots+(i_m-i_{m-1})+(i_m-1)+d_1+\cdots+d_m.
$$ Since $|k-i_1+1|\le|k|+|i_1-1|=k+i_1-1$, the length of $Q$ does
not exceed $k+2(i_m-1)+d_1+\cdots+d_m=|q|+2(i_m-1)<3N$ since
$|q|\le N$, $i_m=\mu(q)\le N$. Now our aim is to construct a
diagram over $\pp$ of area $O(N^2)$ with boundary equation $Q=q$.
We can cancel this equation by $x_0^k$ on the left. That is, we
are going to find the diagram for the following boundary equation:
\be{be} x_0^{-(i_1-1)}x_1^{d_1}x_0^{-(i_2-i_1)}x_1^{d_2}\cdots
x_0^{-(i_m-i_{m-1})}x_1^{d_m}x_0^{i_m-1}=x_{i_1}^{d_1}\cdots
x_{i_m}^{d_m}. \ee

The left-hand side of~(\ref{be}) is a product of two factors. The
second factor is $x_0^s$, where $s=i_m-1$ is non-negative. The
first factor involves $x_0$ in negative powers only. The number of
occurrences of $x_0\iv$ into the first factor equals exactly $s$.
This induces the following definition of a certain class of
diagrams.

Let $\Delta$ be a van Kampen diagram over $\pp$ with boundary
equation of the form $y_1^{-1}uy_2=z$, where $u$, $z$ do not have
occurrences of $x_0^{\pm1}$. Let $y_1$, $y_2$ be decomposed as
$y_j=x_0y_{j1}x_0y_{j2}x_0\ldots y_{j,s-1}x_0y_{js}$, where
$j=1,2$ and no $x_0^{\pm1}$ occur in any of the $y_{ji}$ ($1\le
i\le s$). A diagram with these properties (together with the
decomposition of its contour) will be called {\em balanced}. We
will always draw $u$ on the top, $z$ on the bottom, $y_1$ ($y_2$)
on the left (right).

A very elementary analysis of $\Delta$ allows to conclude that it
contains $s$ subdiagrams that will be called $x_0$-bands. An
$x_0$-band is a diagram with boundary equation of the form
$v^{x_0}=\psi(v)$, where $v$ contains only letters $x_h^{\pm1}$,
$h=1,2,3,4,5$. This $0$-band consists of exactly $|v|$ cells. The
$i$th $x_0$-band in $\Delta$ must connect the $i$th occurrence of
$x_0$ in $y_1$ with the $i$th occurrence of $x_0$ in $y_2$ ($1\le
i\le s$). Although these facts about bands are very easy to prove,
we prefer not to do this in details because all balanced diagrams
we deal with, will already contain the $x_0$-bands with these
properties. So we just include these properties of $x_0$-bands
into the definition of a balanced diagram. The top of the $i$th
$x_0$-band ($1\le i\le s$) will be denoted by
$u_i=u_i(x_1,\ldots,x_5)$. The bottom path of the $i$th band will
thus be $\psi(u_i)=u_i(x_2,\ldots,x_6)$.

Let us define a sequence of subdiagrams $\Theta_0$, $\Theta_1$,
\dots, $\Theta_s$ as follows. Take all the $x_0$-bands off
$\Delta$. The complement of the deleted bands will consist of
$s+1$ subdiagrams that will be enumerated from top to bottom. So
$\Theta_0$ has boundary equation $u=u_1$. For each $1\le i<s$ the
diagram $\Theta_i$ is contained between the $i$th and the
($i+1$)th $x_0$-band. Its boundary equation is
$y_{1i}u_{i+1}=\psi(u_i)y_{2i}$. The boundary equation for
$\Theta_s$ is $y_{1s}z=\psi(u_s)y_{2s}$.

Apply $\psi^i$ to $\Theta_{s-i}$ for each $0\le i\le s$. Clearly,
the bottom label of $\psi^s(\Theta_0)$ equals
$\psi^s(u_1)=\psi^{s-1}(\psi(u_1))$, which is the top label of
$\psi^{s-1}(\Theta_1)$. Analogously, the bottom label of
$\psi^{s-i}(\Theta_i)$ coincides with the top label of
$\psi^{s-i-1}(\Theta_{i+1})$ for all $0\le i<s$. So we can glue
diagrams $\psi^s(\Theta_0)$, $\psi^{s-1}(\Theta_1)$, \dots,
$\psi(\Theta_s)$ together along the paths that have equal labels.
The boundary equation of the result, which will be denoted by
$\Delta'$, has the form $(y_1')^{-1}u'y_2'=z$, where
$u'=\psi^s(u)$,
$y_j'=\psi^{s-1}(y_{j1})\psi^{s-2}(y_{j2})\cdots\psi(y_{j,s-1})y_{js}$
($j=1,2$) with $z$ at the bottom. If we glue $\Delta$ and the
mirror image of $\Delta'$ along the path labelled by $z$, we get
the diagram $\bdelta$ with boundary equation of the form $L=R$,
where \be{coll}
\begin{array}{ccl}
L&=&y_{1s}^{-1}x_0\iv\cdots x_0\iv y_{11}\iv x_0^{-1}
ux_0y_{21}x_0\cdots x_0y_{2,s},\\
R&=&y_{1,s}^{-1}\psi(y_{1,s-1})^{-1}\cdots\psi^{s-1}(y_{11})\iv
\psi^s(u)\psi^{s-1}(y_{21})\cdots\psi(y_{2,s-1})y_{2s}.
\end{array}
\ee

Recall that applying $\psi^k$ means increasing all subscripts by
$k$. The word $R$ in~(\ref{coll}) is the word that can be obtained
as follows. Take the product $y_1\iv uy_2$, delete all occurrences
of $x_0^{\pm1}$ and then increase subscript on each letter by some
number, which is equal to the number of $x_0$-bands that are
contained in $\Delta$ below the corresponding occurrence of the
same letter in $\Delta$. This rule is quite easy to apply.

If $\Delta$ is a balanced diagram, then the we will say that
$\bdelta$ is the result of the {\em collecting process}.
Obviously, the area of $\bdelta$ does not exceed the area of
$\Delta$ multiplied by 2.

Recall that we wanted to find a diagram with boundary
equation~(\ref{be}). Suppose that we found a (balanced) van Kampen
diagram $\Delta$ over $\pp$ of area $O(N^2)$ that has boundary
equation of the form \be{be1}
x_0^{-(i_1-1)}x_1^{d_1}x_0^{-(i_2-i_1)}x_1^{d_2}\cdots
x_0^{-(i_m-i_{m-1})}x_1^{d_m}x_0^{i_m-1}=z, \ee where $z$ is any
word without occurrences of $x_0^{\pm1}$. Let us explain how to
get~(\ref{be}) from there. We have a balanced diagram $\Delta$
with boundary equation $y_1^{-1}uy_2=z$, where $u=x_1^{d_m}$,
$s=i_m-1$, $y_2=x_0^s$, $y_1=x_0^{i_m-i_{m-1}}x_1^{-d_{m-1}}\cdots
x_1^{-d_1}x_0^{i_1-1}$. Collecting process applied to $\Delta$
leads to the diagram $\bdelta$ that has boundary equation $y_1\iv
uy_2=(y_1')\iv u'y_2'$. The word $y_j'$ ($j=1,2$) is obtained from
$y_j$ by deleting all the $x_0$'s and then increasing subscripts
on some letters. Obviously, $y_2'$ is empty whereas each
occurrence of $x_1$ into $y_1'$ is replaced by $x_{1+d}$, where
$d$ is the number of the $x_0$'s that come after the given
occurrence of $x_1$ in $y_1$. The word $u=x_1^{d_m}$ increases all
its subscripts by $s=i_m-1$ so it becomes $u'=x_{i_m}^{d_m}$ (the
subscript 1 must be increased here by
$i_m-1=(i_m-i_{m-1})+\cdots+(i_1-1)$). Analogously,
$x_1^{-d_{m-1}}$ becomes $x_{i_{m-1}}^{-d_{m-1}}$, and so on. As a
result, $(y_1')\iv u'=x_{i_1}^{d_1}\cdots x_{i_m}^{d_m}$. So the
boundary equation of $\bdelta$ is exactly~(\ref{be}). The area of
$\bdelta$ is also $O(N^2)$.

The above reasons show that to prove Lemma~\ref{hor}, it just
remains to find a suitable diagram with boundary
equation~(\ref{be1}). This will be done due to the next Lemma,
which is a crucial point of our proof. \vspace{1ex}

\begin{lm}
\label{crc} There exist positive integer constants $C_1$, $C_2$,
and $D$ such that for any sequence of integers $\alpha_{-n}$,
\dots, $\alpha_{-1}$, $\alpha_0$, $\alpha_1$, \dots, $\alpha_n$,
there exists a word $w=w(x_1,\ldots,x_5)$ and a van Kampen diagram
$\Delta$ over $\pp$ with boundary equation $$
(x_0x_1^{\alpha_{-1}}\cdots x_0x_1^{\alpha_{-n}})\iv
x_1^{\alpha_0} (x_0x_1^{\alpha_1}\cdots
x_0x_1^{\alpha_n})=w(x_1,\ldots,x_5) $$ and the following
conditions hold:

a) $|w|\le S+Dn$,

b) the area of $\Delta$ does not exceed $(C_1S+C_2n)n$, where
$S=\sum_{i=-n}^n|\alpha_i|$.
\end{lm}

To get~(\ref{be1}), it suffices to take $n=i_m-1$. Clearly,
$n<\mu(q)\le N$. The $\alpha_i$'s will be either zero, or they are
equal to $d_m$, $-d_{m-1}$, \dots, $-d_1$. So
$S=d_1+\cdots+d_m\le|q|\le N$. Thus the area of $\Delta$ will be
$O(N^2)$, as desired. So Lemma~\ref{hor} is now finally reduced to
Lemma~\ref{crc}. \vspace{2ex}

{\bf Proof of Lemma~\ref{crc}.}\ We proceed by induction on $n$.
First of all let us construct $\Delta$ for the case $0\le n\le4$.
Consider the following diagram:

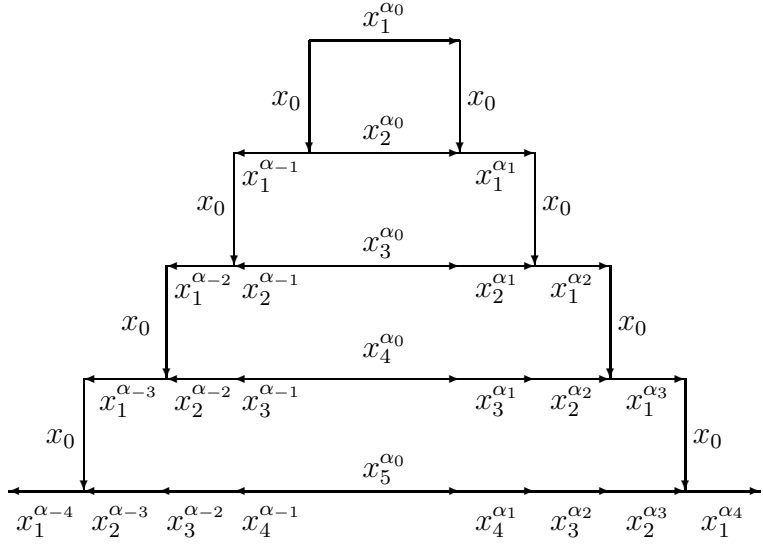
\begin{center}
\unitlength=1mm \special{em:linewidth 0.4pt} \linethickness{0.4pt}
\begin{picture}(102.00,69.00)
\put(42.00,6.00){\vector(1,0){20.00}}
\put(62.00,6.00){\vector(1,0){10.00}}
\put(72.00,6.00){\vector(1,0){10.00}}
\put(82.00,6.00){\vector(1,0){10.00}}
\put(92.00,6.00){\vector(1,0){10.00}}
\put(42.00,6.00){\vector(-1,0){10.00}}
\put(32.00,6.00){\vector(-1,0){10.00}}
\put(22.00,6.00){\vector(-1,0){10.00}}
\put(12.00,6.00){\vector(-1,0){10.00}}
\put(12.00,21.00){\vector(0,-1){15.00}}
\put(42.00,21.00){\vector(1,0){20.00}}
\put(62.00,21.00){\vector(1,0){10.00}}
\put(72.00,21.00){\vector(1,0){10.00}}
\put(82.00,21.00){\vector(1,0){10.00}}
\put(92.00,21.00){\vector(0,-1){15.00}}
\put(42.00,21.00){\vector(-1,0){10.00}}
\put(32.00,21.00){\vector(-1,0){9.00}}
\put(23.00,21.00){\vector(-1,0){11.00}}
\put(42.00,36.00){\vector(1,0){20.00}}
\put(62.00,36.00){\vector(1,0){10.00}}
\put(72.00,36.00){\vector(1,0){10.00}}
\put(82.00,36.00){\vector(0,-1){15.00}}
\put(43.00,36.00){\vector(-1,0){11.00}}
\put(32.00,36.00){\vector(-1,0){9.00}}
\put(23.00,36.00){\vector(0,-1){15.00}}
\put(42.00,51.00){\vector(1,0){20.00}}
\put(62.00,51.00){\vector(1,0){10.00}}
\put(72.00,51.00){\vector(0,-1){15.00}}
\put(42.00,51.00){\vector(-1,0){10.00}}
\put(32.00,51.00){\vector(0,-1){15.00}}
\put(42.00,66.00){\vector(1,0){20.00}}
\put(62.00,66.00){\vector(0,-1){15.00}}
\put(42.00,66.00){\vector(0,-1){15.00}}
\put(52.00,69.00){\makebox(0,0)[cc]{$x_1^{\alpha_0}$}}
\put(52.00,54.00){\makebox(0,0)[cc]{$x_2^{\alpha_0}$}}
\put(52.00,39.00){\makebox(0,0)[cc]{$x_3^{\alpha_0}$}}
\put(52.00,25.00){\makebox(0,0)[cc]{$x_4^{\alpha_0}$}}
\put(52.00,9.00){\makebox(0,0)[cc]{$x_5^{\alpha_0}$}}
\put(65.00,58.00){\makebox(0,0)[cc]{$x_0$}}
\put(75.00,44.00){\makebox(0,0)[cc]{$x_0$}}
\put(85.00,28.00){\makebox(0,0)[cc]{$x_0$}}
\put(95.00,13.00){\makebox(0,0)[cc]{$x_0$}}
\put(39.00,58.00){\makebox(0,0)[cc]{$x_0$}}
\put(29.00,44.00){\makebox(0,0)[cc]{$x_0$}}
\put(19.00,28.00){\makebox(0,0)[cc]{$x_0$}}
\put(9.00,13.00){\makebox(0,0)[cc]{$x_0$}}
\put(7.00,2.00){\makebox(0,0)[cc]{$x_1^{\alpha_{-4}}$}}
\put(17.00,2.00){\makebox(0,0)[cc]{$x_2^{\alpha_{-3}}$}}
\put(27.00,2.00){\makebox(0,0)[cc]{$x_3^{\alpha_{-2}}$}}
\put(37.00,2.00){\makebox(0,0)[cc]{$x_4^{\alpha_{-1}}$}}
\put(67.00,2.00){\makebox(0,0)[cc]{$x_4^{\alpha_1}$}}
\put(77.00,2.00){\makebox(0,0)[cc]{$x_3^{\alpha_2}$}}
\put(87.00,2.00){\makebox(0,0)[cc]{$x_2^{\alpha_3}$}}
\put(97.00,2.00){\makebox(0,0)[cc]{$x_1^{\alpha_4}$}}
\put(37.00,48.00){\makebox(0,0)[cc]{$x_1^{\alpha_{-1}}$}}
\put(67.00,48.00){\makebox(0,0)[cc]{$x_1^{\alpha_1}$}}
\put(37.00,33.00){\makebox(0,0)[cc]{$x_2^{\alpha_{-1}}$}}
\put(67.00,33.00){\makebox(0,0)[cc]{$x_2^{\alpha_1}$}}
\put(37.00,18.00){\makebox(0,0)[cc]{$x_3^{\alpha_{-1}}$}}
\put(67.00,18.00){\makebox(0,0)[cc]{$x_3^{\alpha_1}$}}
\put(28.00,33.00){\makebox(0,0)[cc]{$x_1^{\alpha_{-2}}$}}
\put(77.00,33.00){\makebox(0,0)[cc]{$x_1^{\alpha_2}$}}
\put(28.00,18.00){\makebox(0,0)[cc]{$x_2^{\alpha_{-2}}$}}
\put(77.00,18.00){\makebox(0,0)[cc]{$x_2^{\alpha_2}$}}
\put(18.00,18.00){\makebox(0,0)[cc]{$x_1^{\alpha_{-3}}$}}
\put(87.00,18.00){\makebox(0,0)[cc]{$x_1^{\alpha_3}$}}
\end{picture}
\end{center}

For each $0\le n\le4$, we take the corresponding subdiagram with
bottom label $w(x_1,\ldots,x_5)=(x_n^{\alpha_{-1}}\cdots
x_1^{\alpha_{-n}}) \iv x_{n+1}^{\alpha_0}(x_n^{\alpha_1}\cdots
x_1^{\alpha_n})$. Clearly, $|w|=\sum_{i=-n}^n|\alpha_i|=S$. The
area of $\Delta$ obviously does not exceed $nS$.

From now suppose that $n\ge5$. Let us define integers $k$, $l$,
$m$ in the following way: $k=l=[(n-4)/3]$, $m=n-k-l-4$. Clearly,
$m\ge(n-4)/3>0$. We also have $m\le n-2(n-6)/3-4=n/3$ since
$k=l\ge(n-6)/3$. So the following inequalities hold: \be{fklm}
0\le k,l,m\le\frac{n}3,\quad 2k^2+4l^2+2m^2\le\frac89\,n^2. \ee

Since $0\le k<n$, we can construct a word $u(x_1,\ldots,x_5)$ and
a diagram\footnote{During the proof of this Lemma, by a diagram we
will always mean a van Kampen diagram over $\pp$.} $\Delta_1$ with
boundary equation \be{f1} (x_0x_1^{\alpha_{-1}}\cdots
x_0x_1^{\alpha_{-k}})\iv x_1^{\alpha_0} (x_0x_1^{\alpha_1}\cdots
x_0x_1^{\alpha_k})=u(x_1,\ldots,x_5). \ee Let $S_1$ denote
$\sum_{i=-(k+1)}^{k+1}|\alpha_i|$. (We included the cases
$i=\pm(k+1)$ for convenience.) By the inductive assumption,
\be{est1} |u|\le S_1+Dk,\quad\#\Delta_1\le (C_1S_1+C_2k)k. \ee
Here $\#$ denotes the area of a diagram over $\pp$.

Now let us define 4 new sequences of integers in the following
way:

$$ \beta_i=\alpha_{-(k+i+2)},\ \gamma_i=\alpha_{k+i+2}\ \ (0\le
i\le l+1), $$ and let $\beta_i=\gamma_i=-1$ for all $-l\le i<0$.
Analogously,

$$ \delta_i=\alpha_{-(k+l+i+4)},\ \epsilon_i=\alpha_{k+l+i+4}\ \
(0\le i\le m), $$ and let $\delta_i=\epsilon_i=-1$ for all $-m\le
i<0$.

We will construct 4 diagrams, one for each of these sequences.
They will be denoted by $\Delta_2$, $\Delta_3$, $\Delta_4$,
$\Delta_5$. Let us describe each of these diagrams precisely,
showing its boundary equation and estimating lengths and areas.

The diagram $\Delta_2$ has boundary equation \be{f2}
(x_0x_1^{-1}\cdots x_0x_1^{-1})\iv x_1^{\beta_0}(x_0x_1^{\beta_1}
\cdots x_0x_1^{\beta_l})=u'(x_1,\ldots,x_5). \ee Let $S_2$ denote
$\sum_{i=-l}^{l+1}|\beta_i|=l+\sum_{i=0}^{l+1}|\beta_i|$. (We also
included the case $i=l+1$ for convenience.) By the inductive
assumption, \be{est2} |u'|\le
S_2+Dl,\quad\#\Delta_2\le(C_1S_2+C_2l)l. \ee

The diagram $\Delta_3$ is defined similarly. It has boundary
equation \be{f3} (x_0x_1^{-1}\cdots x_0x_1^{-1})\iv
x_1^{\gamma_0}(x_0x_1^{\gamma_1} \cdots
x_0x_1^{\gamma_l})=u''(x_1,\ldots,x_5). \ee If $S_3$ denotes
$l+\sum_{i=0}^{l+1}|\gamma_i|$, then inequalities \be{est3}
|u''|\le S_3+Dl,\quad\#\Delta_3\le (C_1S_3+C_2l)l \ee hold, as
above.

Finally, $\Delta_4$ and $\Delta_5$ have boundary equations \be{f4}
(x_0x_1^{-1}\cdots x_0x_1^{-1})\iv
x_1^{\delta_0}(x_0x_1^{\delta_1} \cdots
x_0x_1^{\delta_m})=w'(x_1,\ldots,x_5) \ee and \be{f5}
(x_0x_1^{-1}\cdots x_0x_1^{-1})\iv
x_1^{\epsilon_0}(x_0x_1^{\epsilon_1} \ldots
x_0x_1^{\epsilon_m})=w''(x_1,\ldots,x_5), \ee respectively. Let
$S_4=m+\sum_{i=0}^m|\delta_i|$, $S_5=m+\sum_{i=0}^m|\epsilon_i|$.
Then the following estimates hold: \be{est4} |w'|\le
S_4+Dm,\quad\#\Delta_4\le (C_1S_4+C_2m)m, \ee \be{est5} |w''|\le
S_5+Dm,\quad\#\Delta_5\le (C_1S_5+C_2m)m. \ee

Now we shall describe the diagram $\Delta$. First of all, we need
to define the word $w(x_1,\ldots,x_5)$. Let $u_0(x_0,x_1)$,
$u_0'(x_0,x_1)$, and $u_0''(x_0,x_1)$ be the left-hand sides
of~(\ref{f1}), (\ref{f2}), (\ref{f3}), respectively. Let \be{u0}
\bar
u_0(x_1,x_2)=x_1^{-\alpha_{-(k+1)}}u_0(x_1,x_2)x_1^{\alpha_{k+1}},
\ee \be{u00} \bar
u_0'(x_1,x_2)=u_0'(x_1,x_2)x_1^{\beta_{l+1}},\quad \bar
u_0''(x_1,x_2)=u_0''(x_1,x_2)x_1^{\gamma_{l+1}}. \ee By
definition, \be{fw} w(x_1,\ldots,x_5)=w'(x_1,\ldots,x_5)^{-1}\bar
u_0'(x_2,x_3)^{-1} \bar u_0(x_4,x_5)\bar
u_0''(x_2,x_3)w''(x_1,\ldots,x_5). \ee

Each of the diagrams $\Delta_j$ ($1\le j\le5$), $\Delta$ will be
drawn in such a way that words $u$, $w$ (with or without bars
and/or dashes) will be on the bottom. The top will be always
labelled by $x_1$ to a power with zero subscript (directed from
the left to the right). The left (right) side of a diagram will
correspond to the part of its boundary, which involves exponents
with negative (positive) subscript. If $\Gamma$ is a diagram, then
$\Gamma^{-1}$ is a mirror copy of $\Gamma$ with respect to a
horizontal axis symmetry. For a symmetry with respect to a
vertical axis, we will use the notation $-\Gamma$. Recall also
that the operation $\psi$ applied to a word or a diagram increases
all subscripts on the $x_i$'s by 1.

Now we shall describe how the diagram $\Gamma$ is constructed.
First we draw the diagram $\Delta_1$. Its bottom is labelled by
$u(x_1,\ldots,x_5)$. Conjugating this word by $x_0$, gives us an
$x_0$-band that consists of $|u|$ cells. Let us attach this band
to $\Delta_1$. The bottom path of the result is the word
$u(x_2,\ldots,x_6)$. One can then attach $\psi(\Delta_1)^{-1}$ to
it. We get the word $u_0(x_1,x_2)$ on the bottom. Multiplying the
bottom word from both sides by suitable powers of $x_1$, we get
the word $\bar u_0(x_1,x_2)$. Now we conjugate it by $x_0$ and
attach the corresponding $x_0$-band. The word $\bar u_0(x_2,x_3)$
will appear on the bottom of the result.

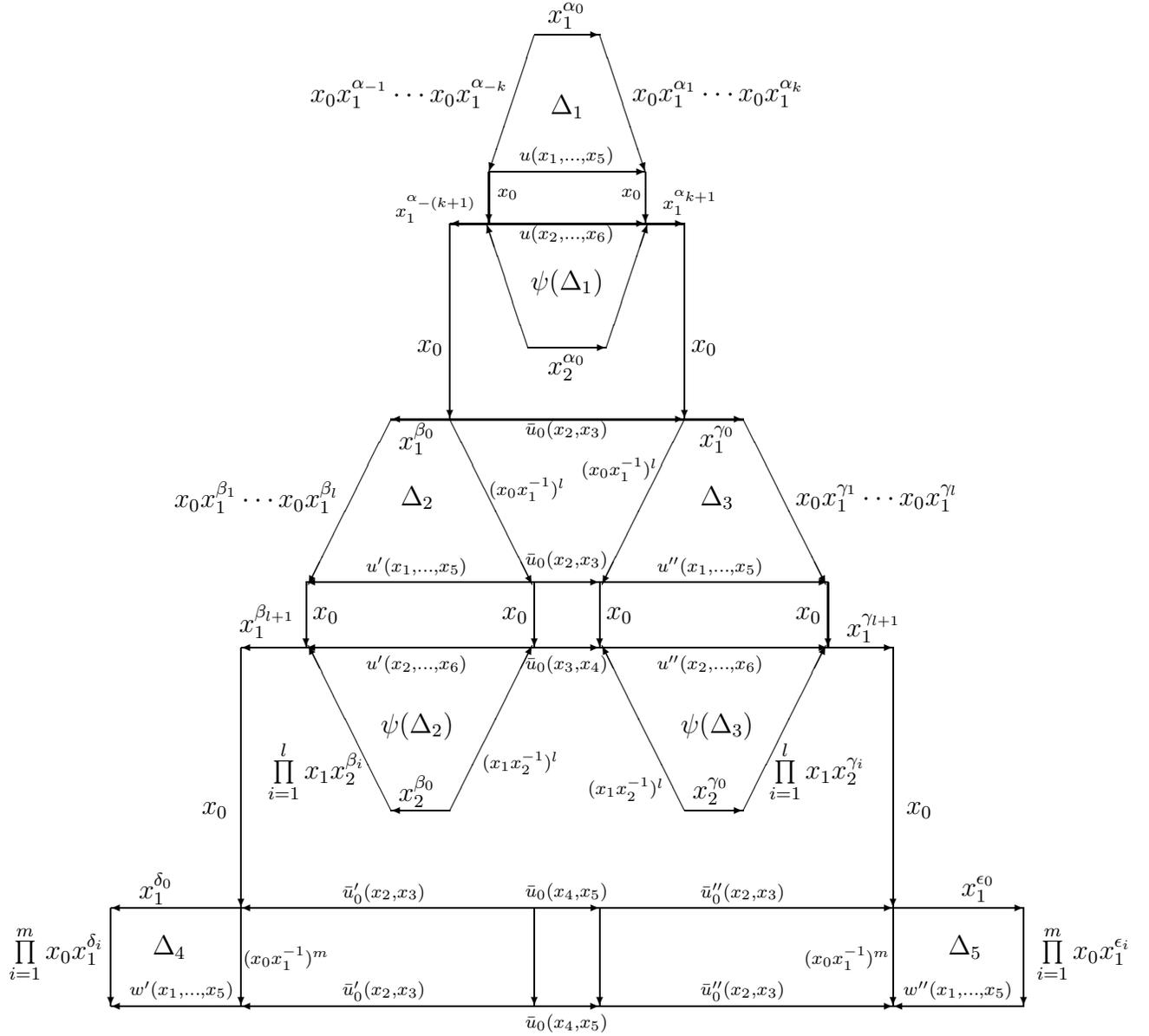
\begin{center}
\unitlength=1.00mm \special{em:linewidth 0.4pt}
\linethickness{0.4pt}
\begin{picture}(148.00,152.00)
\put(26.00,0.00){\vector(-1,0){20.00}}
\put(26.00,0.00){\vector(1,0){100.00}}
\put(146.00,0.00){\vector(-1,0){20.00}}
\put(26.00,15.00){\vector(0,-1){15.00}}
\put(126.00,15.00){\vector(0,-1){15.00}}
\put(26.00,15.00){\vector(1,0){100.00}}
\put(126.00,15.00){\vector(1,0){20.00}}
\put(146.00,15.00){\vector(0,-1){15.00}}
\put(26.00,15.00){\vector(-1,0){20.00}}
\put(6.00,15.00){\vector(0,-1){15.00}}
\put(26.00,55.00){\vector(0,-1){40.00}}
\put(126.00,55.00){\vector(0,-1){40.00}}
\put(36.00,55.00){\vector(-1,0){10.00}}
\put(116.00,55.00){\vector(1,0){10.00}}
\put(71.00,55.00){\vector(1,0){10.00}}
\put(81.00,55.00){\vector(1,0){35.00}}
\put(71.00,55.00){\vector(-1,0){35.00}}
\put(49.00,30.00){\vector(-1,2){12.67}}
\put(58.00,30.00){\vector(1,2){12.67}}
\put(58.00,30.00){\vector(-1,0){9.00}}
\put(94.00,30.00){\vector(-1,2){12.67}}
\put(103.00,30.00){\vector(1,2){12.67}}
\put(94.00,30.00){\vector(1,0){9.00}}
\put(36.00,65.00){\vector(0,-1){10.00}}
\put(71.00,65.00){\vector(0,-1){10.00}}
\put(81.00,65.00){\vector(0,-1){10.00}}
\put(116.00,65.00){\vector(0,-1){10.00}}
\put(71.00,65.00){\vector(1,0){10.00}}
\put(81.00,65.00){\vector(1,0){35.00}}
\put(71.00,65.00){\vector(-1,0){35.00}}
\put(49.00,90.00){\vector(-1,-2){12.67}}
\put(58.00,90.00){\vector(1,-2){12.67}}
\put(58.00,90.00){\vector(-1,0){9.00}}
\put(94.00,90.00){\vector(-1,-2){12.67}}
\put(103.00,90.00){\vector(1,-2){12.67}}
\put(94.00,90.00){\vector(1,0){9.00}}
\put(58.00,90.00){\vector(1,0){36.00}}
\put(58.00,120.00){\vector(0,-1){30.00}}
\put(64.00,120.00){\vector(-1,0){6.00}}
\put(94.00,120.00){\vector(0,-1){30.00}}
\put(88.00,120.00){\vector(1,0){6.00}}
\put(64.00,120.00){\vector(1,0){24.00}}
\put(64.00,128.00){\vector(0,-1){8.00}}
\put(64.00,128.00){\vector(1,0){24.00}}
\put(88.00,128.00){\vector(0,-1){8.00}}
\put(70.00,101.00){\vector(-1,3){6.33}}
\put(82.00,101.00){\vector(1,3){6.33}}
\put(70.00,101.00){\vector(1,0){12.00}}
\put(71.00,149.00){\vector(-1,-3){7.00}}
\put(81.00,149.00){\vector(1,-3){7.00}}
\put(71.00,149.00){\vector(1,0){10.00}}
\put(76.00,152.00){\makebox(0,0)[cc]{$x_1^{\alpha_0}$}}
\put(86.00,140.00){\makebox(0,0)[lc]{$x_0x_1^{\alpha_1}\cdots
x_0x_1^{\alpha_k}$}}
\put(67.00,140.00){\makebox(0,0)[rc]{$x_0x_1^{\alpha_{-1}}\cdots
x_0x_1^{\alpha_{-k}}$}}
\put(76.00,130.00){\makebox(0,0)[cc]{$^{u(x_1,\ldots,x_5)}$}}
\put(76.00,118.00){\makebox(0,0)[cc]{$^{u(x_2,\ldots,x_6)}$}}
\put(76.00,98.00){\makebox(0,0)[cc]{$x_2^{\alpha_0}$}}
\put(97.00,101.00){\makebox(0,0)[cc]{$x_0$}}
\put(55.00,101.00){\makebox(0,0)[cc]{$x_0$}}
\put(64.00,79.00){\makebox(0,0)[lc]{$^{(x_0x_1^{-1})^l}$}}
\put(90.00,82.00){\makebox(0,0)[rc]{$^{(x_0x_1^{-1})^l}$}}
\put(111.00,78.00){\makebox(0,0)[lc]{$x_0x_1^{\gamma_1}\cdots
x_0x_1^{\gamma_l}$}}
\put(41.00,78.00){\makebox(0,0)[rc]{$x_0x_1^{\beta_1}\cdots
x_0x_1^{\beta_l}$}}
\put(99.00,87.00){\makebox(0,0)[cc]{$x_1^{\gamma_0}$}}
\put(53.00,87.00){\makebox(0,0)[cc]{$x_1^{\beta_0}$}}
\put(53.00,67.00){\makebox(0,0)[cc]{$^{u'(x_1,\ldots,x_5)}$}}
\put(98.00,67.00){\makebox(0,0)[cc]{$^{u''(x_1,\ldots,x_5)}$}}
\put(53.00,52.00){\makebox(0,0)[cc]{$^{u'(x_2,\ldots,x_6)}$}}
\put(98.00,52.00){\makebox(0,0)[cc]{$^{u''(x_2,\ldots,x_6)}$}}
\put(53.00,33.00){\makebox(0,0)[cc]{$x_2^{\beta_0}$}}
\put(98.00,33.00){\makebox(0,0)[cc]{$x_2^{\gamma_0}$}}
\put(130.00,30.00){\makebox(0,0)[cc]{$x_0$}}
\put(22.00,30.00){\makebox(0,0)[cc]{$x_0$}}
\put(33.00,7.00){\makebox(0,0)[cc]{$^{(x_0x_1^{-1})^m}$}}
\put(119.00,7.00){\makebox(0,0)[cc]{$^{(x_0x_1^{-1})^m}$}}
\put(139.00,18.00){\makebox(0,0)[cc]{$x_1^{\epsilon_0}$}}
\put(13.00,18.00){\makebox(0,0)[cc]{$x_1^{\delta_0}$}}
\put(17.00,2.00){\makebox(0,0)[cc]{$^{w'(x_1,\ldots,x_5)}$}}
\put(136.00,2.00){\makebox(0,0)[cc]{$^{w''(x_1,\ldots,x_5)}$}}
\put(39.00,60.00){\makebox(0,0)[cc]{$x_0$}}
\put(68.00,60.00){\makebox(0,0)[cc]{$x_0$}}
\put(84.00,60.00){\makebox(0,0)[cc]{$x_0$}}
\put(113.00,60.00){\makebox(0,0)[cc]{$x_0$}}
\put(67.00,124.00){\makebox(0,0)[cc]{$^{x_0}$}}
\put(86.00,124.00){\makebox(0,0)[cc]{$^{x_0}$}}
\put(148.00,8.00){\makebox(0,0)[lc]{$\prod\limits_{i=1}^m
x_0x_1^{\epsilon_i}$}}
\put(123.00,58.00){\makebox(0,0)[cc]{$x_1^{\gamma_{l+1}}$}}
\put(30.00,59.00){\makebox(0,0)[cc]{$x_1^{\beta_{l+1}}$}}
\put(95.00,123.00){\makebox(0,0)[cc]{$^{x_1^{\alpha_{k+1}}}$}}
\put(56.00,122.00){\makebox(0,0)[cc]{$^{x_1^{\alpha_{-(k+1)}}}$}}
\put(76.00,17.00){\makebox(0,0)[cc]{$^{\bar u_0(x_4,x_5)}$}}
\put(71.00,15.00){\vector(0,-1){15.00}}
\put(81.00,15.00){\vector(0,-1){15.00}}
\put(48.00,17.00){\makebox(0,0)[cc]{$^{\bar u_0'(x_2,x_3)}$}}
\put(103.00,17.00){\makebox(0,0)[cc]{$^{\bar u_0''(x_2,x_3)}$}}
\put(76.00,138.00){\makebox(0,0)[cc]{$\Delta_1$}}
\put(76.00,111.00){\makebox(0,0)[cc]{$\psi(\Delta_1)$}}
\put(53.00,78.00){\makebox(0,0)[cc]{$\Delta_2$}}
\put(99.00,78.00){\makebox(0,0)[cc]{$\Delta_3$}}
\put(53.00,43.00){\makebox(0,0)[cc]{$\psi(\Delta_2)$}}
\put(99.00,43.00){\makebox(0,0)[cc]{$\psi(\Delta_3)$}}
\put(71.00,15.00){\vector(1,0){10.00}}
\put(71.00,0.00){\vector(1,0){10.00}}
\put(71.00,15.00){\vector(-1,0){45.00}}
\put(71.00,0.00){\vector(-1,0){45.00}}
\put(76.00,88.00){\makebox(0,0)[cc]{$^{\bar u_0(x_2,x_3)}$}}
\put(76.00,68.00){\makebox(0,0)[cc]{$^{\bar u_0(x_2,x_3)}$}}
\put(76.00,52.00){\makebox(0,0)[cc]{$^{\bar u_0(x_3,x_4)}$}}
\put(5.00,8.00){\makebox(0,0)[rc]{$\prod\limits_{i=1}^m
x_0x_1^{\delta_i}$}} \put(48.00,2.00){\makebox(0,0)[cc]{$^{\bar
u_0'(x_2,x_3)}$}} \put(103.00,2.00){\makebox(0,0)[cc]{$^{\bar
u_0''(x_2,x_3)}$}} \put(76.00,-3.00){\makebox(0,0)[cc]{$^{\bar
u_0(x_4,x_5)}$}}
\put(107.00,36.00){\makebox(0,0)[lc]{$\prod\limits_{i=1}^l
x_1x_2^{\gamma_i}$}}
\put(45.00,36.00){\makebox(0,0)[rc]{$\prod\limits_{i=1}^l
x_1x_2^{\beta_i}$}}
\put(63.00,37.00){\makebox(0,0)[lc]{$^{(x_1x_2^{-1})^l}$}}
\put(91.00,33.00){\makebox(0,0)[rc]{$^{(x_1x_2^{-1})^l}$}}
\put(15.00,9.00){\makebox(0,0)[cc]{$\Delta_4$}}
\put(137.00,9.00){\makebox(0,0)[cc]{$\Delta_5$}}
\end{picture}
\end{center}
\vspace{2ex}

Note that each of the letters $x_2$, $x_3$, $x_4$, $x_5$ commutes
with $x_0x_1\iv$. The diagram of this commutativity consists of
exactly 2 cells of the presentation $\pp$. So we can conjugate
$\bar u_0(x_2,x_3)$ by $(x_0x_1^{-1})^l$ attaching the diagram of
commutativity. Then the same word $\bar u_0(x_2,x_3)$ appears on
the bottom. (The area of the diagram we have attached is exactly
$2l|\bar u_0|$.) Now the word $(x_0x_1^{-1})^l$ appears on the
left and on the right. We are going to attach to it the diagram
$-\Delta_2$ on the left and $\Delta_3$ on the right in a natural
way. After that, the word on the bottom becomes equal to
$u'(x_1,\ldots,x_5)^{-1}\bar u_0(x_2,x_3)u''(x_1,\ldots,x_5)$.
Conjugate this word by $x_0$ and attach an $x_0$-band. The bottom
word will be the image of the previous word under $\psi$, that is,
$u'(x_2,\ldots,x_6)^{-1}\bar u_0(x_3,x_4) u''(x_2,\ldots,x_6)$.
Now we take $\psi(\Delta_2)^{-1}$ and $\psi(\Delta_3)^{-1}$,
attaching them to the subwords $u'$, $u''$, respectively. If we
multiply the result on the left and on the right by suitable
powers of $x_1$, then we get exactly the word $\bar
u_0'(x_1,x_2)^{-1} \bar u_0(x_3,x_4)\bar u''_0(x_1,x_2)$. Adding
an $x_0$-band to the bottom, we increase all subscripts by 1. Then
we get a word in $x_2$, \dots, $x_5$. It commutes with
$(x_0x_1^{-1})^m$ and we add the corresponding diagram of
commutativity. After we attach $-\Delta_4$ on the left and
$\Delta_5$ on the right (gluing along words labelled by
$(x_0x_1^{-1})^m$), we get the desired diagram $\Delta$. Its
bottom path will be labelled by $w(x_1,\ldots,x_5)$, which follows
from~(\ref{fw}). The boundary equation will be exactly as in the
statement, which follows from the definition of the $\beta$'s,
$\gamma$'s, $\delta$'s, and $\epsilon$'s. \vspace{1ex}

Let us estimate the length of $|w|$. From~(\ref{fw}) we see that
$|w|=|w'|+|w''|+|\bar u_0|+|\bar u_0'|+|\bar u_0''|$.
Using~(\ref{est4}), (\ref{est5}), (\ref{u0}), (\ref{u00}), and
taking~(\ref{fklm}) into account, we obtain that $|w|\le
(S_4+Dm)+(S_5+Dm)+(S_1+2k)+(S_2+2l)+(S_3+2l)=
2Dm+(S_1+\cdots+S_5)+2k+4l=2Dm+\sum_{i=-n}^n|\alpha_i|+2l+2m=S+2Dm+(2k+6l+2m)
\le S+2Dn/3+10n/3\le S+Dn$ provided $D\ge10$. (We also refer to
the definition of the $S_j$'s.) So estimate a) holds.

Now we are going to estimate the area of $\Delta$ in the same way.
First of all, we mention that $\Delta$ consists of two copies of
each of the $\Delta_j$, $j=1,2,3$ and one copy of each of the
$\Delta_j$, $j=4,5$. Also $\Delta$ has three $x_0$-bands between
$\Delta_j$ and its mirror copy ($j=1,2,3$). The sum of areas of
these bands clearly equals $|u|+|u'|+|u''|$. The rest of $\Delta$
can be divided into 3 rectangles. The square of each of these
rectangles will be the product of its sides. The first, central,
rectangle will have one side of length $|\bar u_0|$ (width) and
another side of length $2l+2m+3$ (height). The two other
rectangles have width $|\bar u_0'|$, $|\bar u_0''|$, respectively
and the same height $2m+1$. Note that each of the rectangles
indeed consists of small squares $1\times1$, where each square
represents one defining relation of $\pp$. (This explains, by the
way, why we prefer to deal with defining relations of the form
$x_j^{x_i}=x_{j+1}$, where $0<j-i\le5$.) Therefore, the area of
$\Delta$ equals $$
2(\#\Delta_1+\#\Delta_2+\#\Delta_3)+\#\Delta_4+\#\Delta_5+
|u|+|u'|+|u''|+|\bar u_0|(2l+2m+3)+(|\bar u_0'|+|\bar
u_0''|)(2m+1). $$ Note that $|u|\le S_1+Dk$, $|u'|\le S_2+Dl$,
$|u''|\le S_3+Dl$ by (\ref{est1}), (\ref{est2}), (\ref{est3}).
Also recall that $|\bar u_0|= S_1+2k$, $|\bar u_0'|=S_2+2l$,
$|\bar u_0''|=S_3+2l$. We will also use the equality
$S_1+\cdots+S_5=S+2l+2m$, as we did before. Taking into account
the estimates of areas of the $\Delta_j$'s, we finally get
$\#\Delta\le C_1(2S_1k+2S_2l+2S_3l+S_4m+S_5m)+C_2(2k^2+4l^2+2m^2)+
S_1+S_2+S_3+D(k+2l)+(S_1+2k)(2l+2m+3)+(S_2+S_3+4l)(2m+1)$. This
estimate for the area of $\Delta$ can be rewritten as $$
C_2(2k^2+4l^2+2m^2)+S_1(2kC_1+2l+2m+4)+(S_2+S_3)(2lC_1+2m+2)+
(S_4+S_5)C_1m+P, $$ where $P=D(k+2l)+2k(2l+2m+3)+4l(2m+1)\le
Dn+16n^2/9+10n/3$. It follows from the definition of $k$, $l$, $m$
that $l+m+2=n-k-2\le n-(n-6)/3+2\le2n/3$. We also know that $k\le
n/3$ so the coefficient on $S_1$ does not exceed $2nC_1/3+4n/3\le
C_1n$ provided $C_1\ge4$. All coefficients on the other
$\Delta_j$'s satisfy the same inequality. From~(\ref{fklm}) we
also know that the coefficient on $C_2$ does not exceed $8n^2/9$.
Hence the area of $\Delta$ does not exceed $$
8C_2n^2/9+(S_1+\cdots+S_5)C_1n+P=8C_2n^2/9+C_1Sn+
(2l+2m)C_1n+P\le(C_1S+C_2n)n $$ provided $2C_1(l+m)n+P\le
C_2n^2/9$. So it remains to mention that $$
9\,\frac{2(l+m)C_1n+P}{n^2}\le12C_1+\frac{9D+30}n+16\le88 $$ if we
choose $D=10$, $C_1=4$ (recall that $n\ge5$). Thus one can choose
$C_2=88$ to make the proof of part b) complete. \vspace{3ex}

{\bf Remark.}\ We will need a modification of Lemma~\ref{crc} to
prove the Lemma about vertical diagrams. In this case we will need
some flexibility to choose $k$, $l$, $m$. Namely, we need to be
able later to make a choice of $k$ and $m$ from some interval of
length $\xi n$, where $\xi$ is a positive constant. Let us show
how to modify inequalities~(\ref{fklm}) for our purposes. Suppose
that there exist some (small) positive constants $\eta$, $\zeta$
such that inequalities \be{mod} 0\le k,l,m\le(\frac12-\eta)n,\quad
2k^2+4l^2+2m^2\le(1-\zeta)n^2 \ee hold. (In~(\ref{fklm}), we had
$\eta=1/6$, $\zeta=1/9$.) Now for the length of $w$ we will have
$|w|\le S+Dn$ if $2Dm+2k+6l+2m\le Dn$. This obviously holds if $D$
is big enough, say, if $D\ge5/2\eta$. So a) holds.

It is also easy to find out what happens with b). Indeed, to
choose $C_1$ it just suffices to make each coefficient on the
$S_j$'s not to exceed $C_1n$. For the coefficient on $S_1$ we have
$2kC_1+2l+2m+4\le(1-2\eta)nC_1+3n\le nC_1$ if $C_1\ge3/2\eta$. All
the other coefficients satisfy the same inequality. Finally, for
the choice of $C_2$ we need $2C_1(l+m)n+P\le C_2n^2\zeta$. Since
$P$ can be roughly estimated as $4n^2+(3D/2+5)n<(3D/10+5)n^2$, it
suffices to choose some $C_2$ with
$C_2\ge\zeta^{-1}(2C_1+3D/10+5)$. \vspace{2ex}

\section{Vertical Diagrams}
\label{Ve}

The aim of this Section is to prove Lemma~\ref{vert}. We are going
to discuss first how vertical diagrams may look like.

Let $q=x_j$, $p=x_{i_1}x_{i_2}\cdots x_{i_m}$, where $j+m\le n$,
$i_1\le i_2\le\cdots\le i_m$. Suppose that conjugation of $q$ by
$p$ is successful, that is, for any $0\le s<m$, it holds
$j+s>i_{s+1}$. It is easy to construct a diagram over~(\ref{xinf})
of $m$ cells with boundary equation $x_j^p=x_{j+m}$. Indeed, for
any $0\le s<m$, we can conjugate $x_{j+s}$ by $x_{i_{s+1}}$, which
gives $x_{j+s+1}$. All these cells can be concatenated together.
As a result, we get a diagram with boundary equation
$x_j^p=x_{j+m}$. Consider a diagram $\Delta$ over $\pp$ with the
same boundary equation. Let $x_j$ be written on the top and let
$x_{j+m}$ be written on the bottom. We want $\Delta$ to have $m$
subdiagrams $\Gamma_1$, \dots, $\Gamma_m$. For each $1\le s\le m$,
the diagram $\Gamma_s$ is an $x_{i_s}$-band. It ``connects" the
$s$th letter on the left side of $\Delta$ with the $s$th letter on
the right side (these letters have the same label $x_{i_s}$). The
top path of $\Gamma_s$ must have label of the form
$v_s(x_{i_s+1},\ldots,x_{i_s+5})$. The bottom label of $\Gamma_j$
will thus be the image of this word under $\psi$, that is,
$v_s(x_{i_s+2},\ldots,x_{i_s+6})$. It is convenient to define also
diagrams $\Gamma_0$ and $\Gamma_{m+1}$. Each of them consists of
one edge, the top edge and the bottom edge of $\Delta$,
respectively. We also let $i_{m+1}=j+m-1$ by definition.

For each $0\le s\le m$, by $\Theta_s$ we denote the subdiagram of
$\Delta$ that is contained between $\Gamma_s$ and $\Gamma_{s+1}$.
For any $i\ge0$ we will denote by $\bar x_i$ the vector
$(x_i,x_{i+1},x_{i+2},x_{i+3},x_{i+4})$. Our aim is to find a way
to define the words $v_s$ and diagrams $\Theta_s$. The boundary
equation for $\Theta_s$ must be of the form $v_s(\bar x_{i_s+2})
=v_{s+1}(\bar x_{i_{s+1}+1})$ for $1\le s\le m$ (the left-hand
side always means the top label). As for $\Theta_0$, its boundary
equation is $x_j=v_1(\bar x_{i_1+1})$.

Take any diagram $\Xi$ with boundary equation
$x_0^{-n}x_1x_0^n=z$, where $x_1$ is written on the top, $z$ is a
word without occurrences of $x_0^{\pm1}$ written on the bottom.
This diagram must contain $x_0$-bands (recall that this fact can
be proved directly but we just assume this because in our
situation the diagram will automatically have this feature). The
$i$th $x_0$-band connects the $i$th occurrences of $x_0$ into
$x_0^n$ written on the left side and on the right side of the
diagram. Clearly, the top path of the $i$th band can be denoted by
$u_i(\bar x_1)$ for some word $u_i$. The bottom path of the same
band will be $u_i(\bar x_2)$. For any $1\le i<n$, there exists a
subdiagram $T_i$ in $\Xi$ that is contained between the $i$th and
the $(i+1)$th $x_0$-band. Its boundary equation is $u_i(\bar
x_2)=u_{i+1}(\bar x_1)$. We call it (the $i$th) {\em transition\/}
diagram. One can say that $T_i$ converts a word in $x_2$, $x_3$,
\dots into a word in $x_1$, $x_2$, \dots\ . If $T_i$ and $T_{i+1}$
($1\le i\le n-2$) are consecutive transition diagrams, then one
can concatenate $\psi(T_i)$ and $T_{i+1}$ identifying the bottom
path of the first diagram and the top path of the second diagram.
The paths we identify obviously have the same label $u_{i+1}(\bar
x_2)$. The result will be a diagram with boundary equation
$u_i(\bar x_3)=u_{i+2}(\bar x_1)$. This diagram will be denoted as
a product, that is, $\psi(T_i)T_{i+1}$. One can also form multiple
concatenations adjusting them by applying a suitable power of
$\psi$. So we can define products of the form
$$
\psi^k(T_i)\psi^{k-1}(T_{i+1})\cdots\psi^{k-r}(T_{i+r})
$$
for any $i,k,r\ge1$ such that $i+r<n$, $k\ge r$. The boundary equation of
the result will be $u_i(\bar x_{k+2})=u_{i+r+1}(\bar x_{k-r+1})$.
This diagram converts a word in $x_{k+2}$, $x_{k+3}$, \dots into a
word in $x_{k-r+1}$, $x_{k-r+2}$, \dots\,. Thus transition
diagrams can be used to define diagrams $\Theta_s$. We go with the
following definitions. Let

$$ \Theta_0=\psi^{j-2}(T_1)\cdots\psi^{i_1}(T_{j-i_1-1}) $$ (if
$j=1$, then $i_1=0$ so the product is empty, that is, it is a
diagram without cells). The boundary equation for $\Theta_0$ is
\be{t0} u_1(\bar x_j)=u_{j-i_1}(\bar x_{i_1+1}), \ee where the
left-hand side is just $x_j$.

For any $1\le s\le m$ we let $\Theta_s=\psi^{i_s}(T_{j-i_s+s-1})$
if $i_{s+1}=i_s$. We also let $\Theta_s$ be empty in the case
$i_{s+1}=i_s+1$. If $i_{s+1}-i_s\ge2$, then we let $$
\Theta_s=(\psi^{i_{s+1}-1}(T_{j-i_{s+1}+s})\cdots
\psi^{i_s+1}(T_{j-i_s+s-2}))^{-1}. $$ For any $i_{s+1}\ge i_s$ the
diagram $\Theta_s$ will have boundary equation \be{ts}
u_{j-i_s+s-1}(\bar x_{i_s+2})=u_{j-i_{s+1}+s}(\bar x_{i_{s+1}+1}).
\ee

So given a diagram $\Xi$, we construct the diagram $\Delta$ in the
way we have described. (Note that $\Delta$ depends of $p$, $q$,
and $\Xi$.) To find out what is the area of $\Delta$, we need to
define some integers $\kappa_1$, \dots, $\kappa_{n-1}$. Namely, we
let $\kappa_i=k$ whenever $T_i$ or its image under a power of
$\psi$ occurs exactly $k$ times in the $\Theta_s$'s ($0\le s\le
m$). We say that $\kappa_i$ is the {\em intensity\/} of $T_i$. It
is worth noting that the numbers $\kappa_i$ depend of $p$ and $q$
but they are independent on the choice of $\Xi$. Suppose that all
words $u_i$ ($1\le i\le n$) in the diagram $\Xi$ have length
$O(n)$. Then each $\Gamma_s$ ($1\le s\le m$) has area $O(n)$ since
$\#\Gamma_s=|v_s|=O(n)$ (here $v_s$ is some of the $u_i$'s). The
sum of areas of all the $\Gamma_s$'s will be $O(n^2)$ since $m\le
n$. When we take the sum of areas of all the $\Theta_s$ ($0\le
s\le m$), then it will be the sum of areas of the $T_i$'s, where
each $T_i$ occurs exactly $\kappa_i$ times. So

\be{aread} \#\Delta=\sum_{i=1}^{n-1}\kappa_i\#T_i+O(n^2). \ee

Let us estimate $\sum_{i=1}^{n-1}\kappa_i$. It is the total number
of factors of the form $T_i$ ($1\le i<n$) for all the
$\Theta_s$'s. By definition, $\Theta_0$ has $j-i_1+1$ factors and
each $\Theta_s$ for $1\le s\le m$ has $|i_{s+1}-i_s-1|$ factors.
So

\be{sumkap}
\sum_{i=1}^{n-1}\kappa_i=j-i_1+1+\sum_{i=1}^{n-1}|i_{s+1}-i_s-1|.
\ee

Since $|i_{s+1}-i_s-1|\le i_{s+1}-i_s+1$, we can
estimate~(\ref{sumkap}) as $j+n+i_{m+1}<2j+m+n\le3n$. So

\be{skon} \sum_{i=1}^{n-1}\kappa_i=O(n). \ee

The diagram $\Xi$ can be easily constructed by Lemma~\ref{crc}.
However, it is clear that diagrams $T_i$ may sometimes have
relatively big area. Although the area of the whole diagram $\Xi$
will be $O(n^2)$, for some values of $i$ the area of $T_i$ may
exceed $cn^2$ for some constant $c>0$. Also we have no guarantee
that $\kappa_i$ will be small. It may happen that $\kappa_i>dn$
for some constant $d>0$. So we never get the estimate $O(n^2)$ if
we choose $\Xi$ arbitrarily. To avoid this difficulty, we need to
construct $\Xi$ in such a way that if $\kappa_i$ is ``big" then
the area of $T_i$ is ``small". This is always possible by the
following Lemma.

\begin{lm}
\label{intens} There exist constants $C,C_0,D>0$ such that for any
integers $\alpha_i$, where $|\alpha_i|\le1$ for all $-n\le i\le n$
and for any non-negative integers $\kappa_i$ $(1\le i\le n)$,
there exists a van Kampen diagram $\Delta'$ over $\pp$ with
boundary equation of the form $$ L\iv x_1^{\alpha_0}R=w(\bar x_1),
$$ where $L=x_0x_1^{\alpha_{-1}}\cdots x_0x_1^{\alpha_{-n}}$,
$R=x_0x_1^{\alpha_1}\cdots x_0x_1^{\alpha_n}$ with the following
properties:

a) $\Delta'$ contains $x_0$-bands $\Delta_i'$ $(1\le i\le n)$,
where $\#\Delta_i'\le Dn$ for all $1\le i\le n$. The $i$th
$x_0$-band connects the $i$th occurrence of $x_0$ into $L$ with
the $i$th occurrence of $x_0$ into $R$. Besides, $|w|\le Dn$.

b) $\#\Delta'\le C_0n^2$.

c) Let $\Phi_i$ $(1\le i\le n)$ denote the subdiagram in $\Delta'$
that is contained between $\Delta_i'$ and $\Delta_{i+1}'$, where
$\Delta_{n+1}'$ denotes the bottom path of $\Delta'$ labelled by
$w(\bar x_1)$. Then $$ \sum_{i=1}^n\kappa_i\#\Phi_i\le
C\sum_{i=1}^n\kappa_i\cdot n. $$
\end{lm}

First of all, let us show that Lemma~\ref{intens} implies what we
need. Take $L=x_0(x_0x_1^{-1})^{n-1}$, $R=x_0^n$, $\alpha_0=0$ and
consider the diagram $\Delta'$ from Lemma~\ref{intens}. (The
integers $\kappa_i$ for $1\le i<n$ are defined by words $p$, $q$
from the statement of Lemma~\ref{vert}; for $\kappa_n$ we can take
$0$.) The diagram $\Delta'$ has boundary equation $Lw=R$. It is
easy to see that $L$ is a smooth word of rank 1. Conjugating $x_1$
by $L$, we get $x_2$. There exists a natural van Kampen diagram
over $\pp$ with boundary equation $x_1^L=x_2$. It has $|L|=2n-1$
cells. We can attach two copies of $\Delta'$ to it. This gives a
diagram with boundary equation $x_1^{x_0^n}=x_2^w$. Denote this
diagram by $\Xi$. As above, we find $x_0$-bands in $\Xi$ and
define transition diagrams $T_i$ ($1\le i<n$). It is now easy to
describe how they look like. Let $w_i(\bar x_1)$ be the top label
of $\Delta_i'$ for each $1\le i\le n$. Then $\Phi_i$ ($1\le i<n$)
has boundary equation $x_1w_i(\bar x_2)=w_{i+1}(\bar x_1)$. It is
easy to see that $T_1$ consists of just one edge labelled by
$x_2$. For each $1<i<n$, $T_i$ has boundary equation $w_i(\bar
x_2)^{-1}x_3w_i(\bar x_2)= w_{i+1}(\bar x_1)^{-1}x_2w_{i+1}(\bar
x_1)$. This diagram consists of one cell of the form $x_3=x_1\iv
x_2x_1$ and two copies of $\Phi_i$ from both sides:

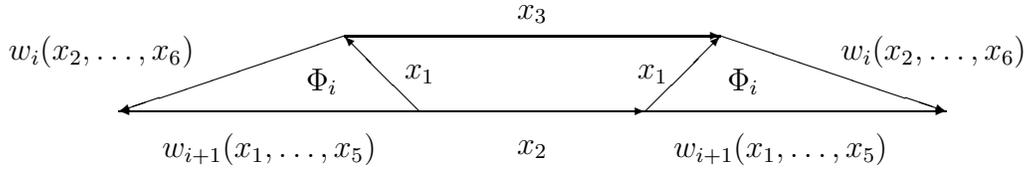
\begin{center}
\unitlength=1mm \special{em:linewidth 0.4pt} \linethickness{0.4pt}
\begin{picture}(112.00,20.00)
\put(32.00,17.00){\vector(1,0){50.00}}
\put(42.00,7.00){\vector(1,0){30.00}}
\put(72.00,7.00){\vector(1,1){10.00}}
\put(42.00,7.00){\vector(-1,1){10.00}}
\put(82.00,17.00){\vector(3,-1){30.00}}
\put(72.00,7.00){\vector(1,0){40.00}}
\put(42.00,7.00){\vector(-1,0){40.00}}
\put(32.00,17.00){\vector(-3,-1){30.00}}
\put(57.00,20.00){\makebox(0,0)[cc]{$x_3$}}
\put(57.00,2.00){\makebox(0,0)[cc]{$x_2$}}
\put(90.00,2.00){\makebox(0,0)[cc]{$w_{i+1}(x_1,\ldots,x_5)$}}
\put(22.00,2.00){\makebox(0,0)[cc]{$w_{i+1}(x_1,\ldots,x_5)$}}
\put(98.00,15.00){\makebox(0,0)[lc]{$w_i(x_2,\ldots,x_6)$}}
\put(12.00,15.00){\makebox(0,0)[rc]{$w_i(x_2,\ldots,x_6)$}}
\put(42.00,12.00){\makebox(0,0)[cc]{$x_1$}}
\put(73.00,12.00){\makebox(0,0)[cc]{$x_1$}}
\put(85.00,11.00){\makebox(0,0)[cc]{$\Phi_i$}}
\put(29.00,11.00){\makebox(0,0)[cc]{$\Phi_i$}}
\end{picture}
\end{center}

So by~(\ref{aread}) the area of $\Delta$ (constructed by $p$, $q$,
$\Xi$) does not exceed $$
\sum_{i=1}^{n-1}\kappa_i(2\#\Phi_i+1)+O(n^2)\le\sum_{i=1}^{n-1}
\kappa_i\cdot(2Cn+1)+O(n^2)=O(n^2) $$ by part c) of
Lemma~\ref{intens} and~(\ref{skon}).

To prove Lemma~\ref{vert}, it suffices to define the smooth path
$t$ from the statement and establish some of its properties. Each
diagram $T_i$ can be naturally decomposed into three parts. In our
case, $T_i$ for $1<i<n$ consists of the central part which has one
cell and two other parts, the left part and the right part, which
are mirror copies of each other (and are equal to $\Phi_i$). The
diagram $T_1$ consists of its central part only (the left part and
the right part of it are empty). When we concatenate some images
of the $T_i$'s under a power of $\psi$, the central parts of the
factors are concatenated too. So each $\Theta_s$ ($0\le s\le m$)
can be naturally decomposed into three parts as well. By $t_s$
($0\le s\le m$) we denote the word written on the left and the
right side of the central part of $\Theta_s$ (these sides must
have the same label by definition).

For any diagram that is decomposed naturally into the three parts,
we can define its {\em thickness\/} as the length of the word
written on the left and the right side of the central part. All
our transition diagrams $T_i$ ($1\le i<n$) have thickness at most
one. So the thickness of $\Theta_s$ does not exceed the sum of its
factors for any $0\le s\le m$. This implies that
$\sum_{s=0}^m|t_s|$ does not exceed $\sum_{i=1}^{n-1} \kappa_i$,
which is $O(n)$ by~(\ref{skon}). To decompose $\Delta$ into three
parts, one need to find this decomposition for each of the bands
$\Gamma_s$ ($1\le s\le m$). Each word $v_s$ ($1\le s\le m$) is one
of the $u_i$'s so it has a ``central" letter. So the band
$\Gamma_s$ must have a ``central" cell. We now can define the word
$t$ as $$ t=t_0x_{i_1}t_1x_{i_2}\cdots x_{i_m}t_m. $$ The diagram
$\Delta$ decomposes into three parts, where the central part $B$
has boundary equation $t^{-1}x_jt=x_{j+m}$ and the two other parts
are mirror copies of each other and they have boundary equations
$q=t$. Since $|t|=\sum_{s=0}^m|t_s|+m=O(n)+m=O(n)$, the only thing
we need to check is that $t$ is a smooth word of rank $j$, height
$m$. (It has been already shown that the area of $\Delta$ is
$O(n^2)$ so the same is true for the area of the subdiagram of
$\Delta$ with boundary equation $q=t$.)

It is clear that $B$ consists of exactly $|t|$ cells. The boundary
equation of each of them has the form $x_k^{x_i}=x_{k+1}$ for some
$k>i$. The letters $x_i$ here belong to the left and the right
side of $B$. So what we really have to prove is the alternative
$k=i+1$ or $k=i+2$. Each letter of $t$ is associated with exactly
one cell so we need to establish this alternative for each letter
of $t$. Each cell that corresponds to a letter from any of the
$t_s$ ($0\le s\le m$) obviously has this property because it is an
image of a central cell of some transition diagram under a power
of $\psi$. But we know that any defining relation for a central
cell of a transition diagram is $x_2^{x_1}=x_3$. So if we apply a
power of $\psi$, then the result will be of the form
$x_k^{x_i}=x_{k+1}$, where $k=i+1$. So we need to check out only
the cells of $B$ that correspond to letters $x_{i_s}$ ($1\le s\le
m$).

Recall that by $u_i(x_1,\ldots,x_5)$ we denoted the top path of
the $i$th $x_0$-band of $\Xi$. The central letter of $u_1$ is
$x_1$; for all $1<i\le n$ the central letter of
$u_i(x_1,\ldots,x_5)$ is $x_2$. Now we can look at the boundary
equations of the diagrams $\Theta_s$ ($0\le s\le m$). According
to~(\ref{t0}), the central letter of the bottom path of $\Theta_0$
is the central letter of the word $u_{j-i_1}(\bar x_{i_1+1})$.
Thus it is $x_{i_1+1}$ if $j-i_1=1$ and it is $x_{i_1+2}$
otherwise. The defining relation that corresponds to the letter
$x_{i_1}$ has the form $x_k^{x_i}=x_{k+1}$, where $i=i_1$ and
$x_k$ is the central letter of the bottom path of $\Theta_0$. This
shows that $k=i+1$ if $j=i_1+1$ and $k=i+2$ if $j-i_1\ge2$.

Now consider the boundary equation of $\Theta_s$, where $1\le
s<m$. According to~(\ref{ts}), the bottom path of $\Theta_s$ is
$u_{j-i_{s+1}+s}(\bar x_{i_{s+1}+1})$. We have $j+s>i_{s+1}$ by
the conditions of Lemma~\ref{vert}. If $j-i_{s+1}+s=1$, then the
central letter we deal with is $x_{i_{s+1}+1}$, otherwise it is
$x_{i_{s+1}+2}$. So the cell that corresponds to $x_{i_{s+1}}$ has
the form $x_k^{x_i}=x_{k+1}$, where $i=i_{s+1}$, $k=i_{s+1}+1$ or
$k=i_{s+1}+2$. This shows that our fact is true for letters
$x_{i_2}$, \dots, $x_{i_m}$.

So $t$ is smooth of rank $j$. Its height is equal to $m$ since
$x_j^t=x_{j+m}$. This completes the proof of Lemma~\ref{vert}.
\vspace{2ex}

Now we are going to prove Lemma~\ref{intens}. It is clear that for
any $n$ there exists at least one diagram with boundary equation
of the Lemma. Let $n_0$ be some positive integer. We can assume
that for any $n<n_0$ all conditions a) -- c) hold. Indeed, the
number of our diagrams under the restriction $n<n_0$ is finite
because all the $\alpha_i$'s are either $0$ or $\pm1$. So if the
constants are chosen big enough, then all inequalities in a) -- c)
hold.

Let $n\ge n_0$. We proceed by induction on $n$. Now we need to
define $k$, $l$, $m$ such that $k+l+m+4=n$ and for some positive
$\eta$, $\zeta$ conditions~(\ref{mod}) hold. We need some freedom
of choice of $k$ and $m$. In our case it is very important to
choose $k$ and $m$ from an interval of length at least $cn$ for
some positive constant $c$. Let $k$, $m$ satisfy the following
inequalities: $$ n'=\left[\frac{n}3\right]\le k,m\le
\left[\frac{n}3\right]+\left[\frac{n}{30}\right]=n''. $$ Thus each
of the $k$, $m$ can take at least $[n/30]+1>n/30$ integer values.
Then we define $l=n-k-m-4$. We need $l\ge0$. Clearly, $l\ge
n-2(n/3+n/30)-4=4(n/15-1)\ge0$ if we take $n_0=15$. On the other
hand, $l\le n-2[n/3]-4\le n-2(n-2)/3-4<n/3$. So each of the $k$,
$l$, $m$ does not exceed $11n/30$ and we can let $\eta=2/15$.
Finally, $2k^2+4l^2+2m^2<4n^2((11/30)^2+1/9)=221n^2/225$ so let
$\zeta=4/225$.

The remark after the proof of Lemma~\ref{crc} shows that for any
$n\ge n_0$ we can choose any $k$ and $m$ such that $[n/3]\le
k,m\le[n/3]+[n/30]$. Thus all estimates of Lemma~\ref{crc} hold.
In particular, we have that the area of $\Delta'$ will be
$O(n^2)$. (Note that $S=\sum_{i=-n}^n|\alpha_i|\le2n+1$ in our
case.) Therefore, part b) holds.

It is easy to prove part a). Take any $x_0$-band $\Delta_i'$. If
$i\le k$, then it is contained in $\Delta_1$. By the inductive
assumption, its area does not exceed $Dk<Dn$. The $(k+1)$th band
has area $|u|\le Dk<Dn$. The ($k+2$)nd band has area $|\bar
u_0|\le4k+3<Dn$ if $D$ is big enough. The next $l$ bands consist
of two bands that are contained in $\Delta_2$, $\Delta_3$,
respectively and the third part that has area $|\bar u_0|$. By the
assumption, the total area will be at most $2Dl+4k+3<2Dn/3+2n\le
Dn$ if $D\ge6$. The ($k+l+3$)rd band satisfies the same condition
since words $u'$, $u''$ also have length at most $Dl$ each. The
next band has area $|\bar u_0'|+|\bar u_0''|+|\bar
u_0|\le4k+8l+7<5n\le Dn$ if $D\ge5$. The area of each of the last
$m$ bands and the length of $w$ can be estimated from above as
$2Dm+4k+8l+7<2Dm+5n\le(11D/15+5)n<Dn$ if $D\ge20$.

It remains to prove the most difficult part c). Note that we have at least two
values of $i$ such that $\Phi_i$ may have big area. For $i=k+1$ the area of
$\Phi_i$ equals $\#\Delta_1\le C_0k^2$. Also for $i=k+l+3$ the area of $\Phi_i$
equals $\#\Delta_2+\#\Delta_3\le2C_0l^2$. So $k$, $l$, $m$ should be chosen in
such a way that $\kappa_i$ would be small for these values of $i$. We did not
specify yet how $k$, $m$ must be chosen. Now we do that. Let $$
\kappa_{k+1}=\min\left\{\,\kappa_{i+1}\mid n'\le i\le n'' \,\right\}, $$ $$
\kappa_{n-m-1}=\min\left\{\,\kappa_{n-i-1}\mid n'\le i\le n''\,\right\}. $$
Since each of the $k$, $m$ may have at least $n/30$ values, it follows from
these definitions that
$$
\kappa_{k+1}\le30\cdot\frac{\sum\limits_{i=n'}^{n''}\kappa_i}n\le
\frac{30\kappa}n,
$$
where $\kappa=\sum_{i=1}^n\kappa_i$. Analogously,
$$
\kappa_{k+l+3}=\kappa_{n-m-1}\le\frac{30\kappa}n.
$$
So the most ``dangerous" part of the sum $\sum_{i=1}^n\kappa_i\#\Phi_i$
can be estimated as follows:
$$
\kappa_{k+1}\#\Phi_{k+1}+\kappa_{k+l+3}\#\Phi_{k+l+3}\le30C_0
\kappa\cdot\frac{k^2+2l^2}n<15C_0\kappa n.
$$

By the inductive assumption, $\sum_{i=1}^k\kappa_i\#\Phi_i\le
Ck\sum_{i=1}^k\kappa_i$. Note that $\Phi_i$ has zero area for
$i=k+2$ and $i=k+l+4$. For each $k+3\le i\le k+l+2$ the diagram
$\Phi_i$ consists of three parts: $\Phi_i'$, which is a subdiagram
of $\Delta_2$, $\Phi_i''$, which is a subdiagram of $\Delta_3$,
and the third part that has area $|\bar u_0|$. By the inductive
assumption applied to $\Delta_2$ and $\Delta_3$, one can assume
that $\sum\limits_{i=k+3}^{k+l+2}\kappa_i\#\Phi_i'\le
Cl\sum\limits_{i=k+3}^{k+l+2}\kappa_i$ and
$\sum\limits_{i=k+3}^{k+l+2}\kappa_i\#\Phi_i''\le
Cl\sum\limits_{i=k+3}^{k+l+2}\kappa_i$. Since $|\bar u_0|\le Dn$,
we have $\sum\limits_{i=k+3}^{k+l+2}\kappa_i\#\Phi_i
\le(2Cl+Dn)\sum\limits_{i=k+3}^{k+l+2}\kappa_i\le(C(1-2\eta)+D)n
\sum\limits_{i=k+3}^{k+l+2}\kappa_i$.

Analogously, for each $k+l+5\le i\le n$ we can naturally subdivide
$\Phi_i$ into three parts. The third part will have area $|\bar
u_0'|+|\bar u_0''|+|\bar u_0|\le Dn$. So we get the inequality
$\sum_{i=k+l+5}^n\kappa_i\#\Phi_i\le(2Cm+Dn)\sum_{i=k+l+5}^n
\kappa_i\le(C(1-2\eta)+D)n\sum_{i=k+l+5}^n\kappa_i$. Taking the
sum of the estimates for each of the cases: $1\le i\le k$, $k+3\le
i\le k+l+2$, $k+l+5\le i\le n$, and $i\in\{\,k+1,k+l+3\,\}$, we
obtain $$
\sum_{i=1}^n\kappa_i\#\Phi_i\le(C(1-2\eta)+D+15C_0)\kappa n\le
C\kappa n $$ if $C\ge(D+15C_0)/2\eta$. This proves part c).
\vspace{1ex}

Lemma~\ref{intens} is proved. Now the main result follows.
\vspace{2ex}

Note that we also got polynomial isoperimetric inequalities for
the two other R.\,Thomp\-son groups, $T$, and $V$ in~\cite{Gu00}
(both groups are finitely presented and simple). Results of the
present paper imply some improvements to the estimates
from~\cite{Gu00}. Details will appear elsewhere.

Note that the techniques of this paper may probably be applied to
other Thompson-like groups. It looks like similar results can be
obtained without any essential changes for groups $F_r$, where
$r\ge3$ (see~\cite{Bro}). It is also interesting to look from this
point of view to some diagram groups (see~\cite{GbS}), especially
for the diagram groups over semigroup presentation $\la x\mid
x^3=x^2\ra$.

\end{document}